\documentclass[reqno]{amsart}
\usepackage{enumerate, bbm}
\usepackage{amssymb,url,color, booktabs}
\usepackage{mathrsfs}
\usepackage[nobysame,abbrev]{amsrefs}
\BibSpec{article}{%
+{}{\PrintAuthors} {author}
+{,}{ \textrm} {title}
+{.}{ \textit} {journal}
+{,}{ \textbf} {volume}
+{}{ \parenthesize} {date}
+{,}{ } {pages}
+{.}{ arXiv:} {eprint}
+{.}{} {transition}
}
\BibSpec{book}{%
+{}{\PrintAuthors} {author}
+{,}{ \textit} {title}
+{.}{ \textrm} {series} 
+{,}{ Vol.} {volume} 
+{.}{ } {publisher}
+{,}{ } {date}
+{,}{ } {pages}
+{.}{} {transition}
}
\usepackage{color}
\usepackage[colorlinks=true]{hyperref}
\hypersetup{
    linkcolor=blue,          
    citecolor=red,        
    filecolor=blue,      
    urlcolor=cyan
}
\definecolor{MyDarkBlue}{cmyk}{0.8,0.3,0.8,0.4}
\definecolor{yellow}{rgb}{0.99,0.99,0.70}
\definecolor{white}{rgb}{1.0,1.0,1.0}
\definecolor{black}{rgb}{0.00,0.00,0.00}
\definecolor{backgroundcolor}{RGB}{199,238,206}

\numberwithin{equation}{section}

\newtheorem{theorem}{Theorem}[section]
\newtheorem{lemma}[theorem]{Lemma}
\newtheorem{remark}[theorem]{Remark}
\newtheorem{definition}[theorem]{Definition}
\newtheorem{proposition}[theorem]{Proposition}
\newtheorem{Examples}[theorem]{Example}
\newtheorem{corollary}[theorem]{Corollary}

\def\eps{\varepsilon}
\def\e{\mathrm{e}}
\def\supp{\mathrm{supp}}

\def\Law{{\mathord{{\rm Law}}}}
\def\dif{{\mathord{{\rm d}}}}

\def\bbone{{\boldsymbol{1}}}
\def\bb2{{\boldsymbol{2}}}
\def\no{\nonumber}
\def\={&\!\!=\!\!&}

\def\bx{{\mathbf{x}}}

\def\bC{{\mathbf C}}

\def\bP{{\mathbf P}}

\def\b1{{\mathbbm 1}}

\def\cB{{\mathcal B}}

\def\cF{{\mathcal F}}

\def\cI{{\mathcal I}}

\def\cR{{\mathcal R}}

\def\mB{{\mathbb B}}

\def\mE{{\mathbb E}}

\def\mL{{\mathbb L}}

\def\mN{{\mathbb N}}

\def\mP{{\mathbb P}}

\def\mR{{\mathbb R}}

\def\sF{{\mathscr F}}

\def\sI{{\mathscr I}}

\def\sP{{\mathscr P}}

\def\sS{{\mathscr S}}

\def\geq{\geqslant}
\def\leq{\leqslant}
\def\<{{\langle}}
\def\>{{\rangle}}
\def\({{\big(}}
\def\){{\big)}}

\def\a{\alpha}
\def\b{\beta}
\def\de{\delta}

\def\l{\lambda}
\def\om{\omega}
\def\Om{\Omega}
\def\s{\sigma}

\def\ff{\frac}
\def\nn{\nabla}

\def\p{\partial}
\def\div{\mathord{{\rm div}}}

\def\bt{\begin{theorem}}
\def\et{\end{theorem}}
\def\bl{\begin{lemma}}
\def\el{\end{lemma}}
\def\br{\begin{remark}}
\def\er{\end{remark}}
\def\bpf{\begin{proof}}
\def\epf{\end{proof}}
\def\bx{\begin{Examples}}
\def\ex{\end{Examples}}
\def\bd{\begin{definition}}
\def\ed{\end{definition}}
\def\bp{\begin{proposition}}
\def\ep{\end{proposition}}
\def\bc{\begin{corollary}}
\def\ec{\end{corollary}}

\def\wt{\widetilde}

\allowdisplaybreaks

\begin{document}

\title{Well-Posedness for SDEs with Logarithmical Critical Distributional Drifts}

\author{Zikai Chen, Zimo Hao and Xicheng Zhang}

\date{\today}


\thanks{This work is supported by National Key R\&D program of China (No. 2023YFA1010103) and NNSFC grant of China (No. 12131019) and the DFG through the CRC 1283 ``Taming uncertainty and profiting from randomness and low regularity in analysis, stochastics and their applications''. Z. Chen was also supported by the China Scholarship Council (Grant No. 202506270058), and Z. Hao by the NNSFC (Grant No. 12601259).}

\begin{abstract}
We study the stochastic differential equation
$$\dif X_t=b(t,X_t)\dif t+\sqrt{2}\dif W_t$$
on $\mR^d$, where $b$ is a time-dependent, divergence-free distributional drift of critical H\"older--Besov regularity $-1$, strengthened by an iterated-logarithmic correction. For every initial probability law, we construct a weak solution by smooth approximation and realize the singular drift as an additive functional. The main analytic ingredient is the Schauder estimate with a logarithmic smallness factor. Combined with uniform logarithmic Krylov estimates and a stochastic substitution formula for distributional test functions, this estimate allows us to apply a Zvonkin transformation and prove uniqueness in law among weak solutions satisfying the corresponding Krylov bounds. For solutions starting from deterministic points, we further show that their time-marginal distributions admit densities satisfying two-sided Aronson-type Gaussian estimates. 
\end{abstract}

\maketitle

\section{Introduction}

Let $T>0$ and $d\geq2$. We consider the following stochastic differential equation on $\mR^d$ driven by a standard $d$-dimensional Brownian motion $W$: for $t\in[0,T]$
\begin{align}\label{SDE}
    \dif X_t=b(t,X_t)\dif t+\sqrt{2}\dif W_t,\ X_0\sim\mu\in\sP(\mR^d),
\end{align}
where $\sP(\mR^d)$ is the space of all probability measures on $\mR^d$. The drift $b$ is allowed to be a time-dependent distribution of critical order $-1$. Thus the expression $b(t,X_t)$ is not meaningful pointwise. Following the approximation approach in \cite{HZ25}, we define the drift term through mollification. Let $\phi$ be a compactly supported smooth probability density and set $\phi_n(x):=n^d\phi(nx)$. For a distributional drift $b$, define
\begin{align}
    b_n(t,x):=(b(t,\cdot)*\phi_n)(x),\no
\end{align}
where $*$ denotes convolution in the spatial variable in the distributional sense.

We first specify the notion of weak solution used in this paper. Since $b(t,X_t)$ is not defined pointwise, the usual integral formulation of the SDE has to be replaced by an approximation-based one.

\begin{definition}[Weak solution]\label{weaksol}
Let $\mathfrak F:=(\Omega,\sF,(\sF_t)_{t\geq0},\bP)$ be a stochastic basis, and let $(X,W)$ be a pair of $\mR^d$-valued continuous $\sF_t$-adapted processes on $\mathfrak F$. We call $(\mathfrak F,X,W)$ a weak solution of SDE \eqref{SDE} with initial distribution $\mu\in\sP(\mR^d)$ if $W$ is an $\sF_t$-Brownian motion, $\bP\circ X_0^{-1}=\mu$, and
\begin{align}
    X_t=X_0+A_t^b+\sqrt{2}W_t,\  t\in[0,T],\ a.s.,\no
\end{align}
where $A_t^b:=\lim_{n\to\infty}\int_0^t b_n(s,X_s)\dif s$ exists in the $L^2$-sense.
\end{definition}

\begin{remark}
It should be noted that the definition of a weak solution in Definition \ref{weaksol} relies on the choice of mollifiers $(\phi_n)_{n\geq1}$. However,  the logarithmic Krylov estimates obtained below imply that the additive functional $A^b$ is canonically determined by $b$ along solutions in the Krylov class.
\end{remark}

\subsection{Main result}
The main result is as follows. The logarithmically corrected H\"older--Besov spaces $\bC^{\a,*}_{m,\eps}$ and $\bC^{\a,**}_{m,\eps}$ are introduced in Section \ref{sec21} below, with $\ell_m$ defined in \eqref{ellm}.

\bt\label{main-result}
Assume that $\eps>0$, $m\in\mN$, $b\in\mL_T^\infty\bC^{-1,*}_{m,\eps}$ and $\div b=0$. Then for every $\mu\in\sP(\mR^d)$, there exists a weak solution $(\mathfrak F,X,W)$ to SDE \eqref{SDE} in the sense of Definition \ref{weaksol}, which is also unique in the class such that the logarithmic Krylov holds: for every $p>2$ and every $f\in\mL_T^1C_b\cap\mL_T^\infty\bC^{-1,**}_{m,\eps}$, there exists a constant $C=C(p,T,d,m,\eps,\|b\|_{\mL_T^\infty\bC^{-1,*}_{m,\eps}})>0$ such that for any $0\le t_0<t_1\le T$,
\begin{align}
    \Big\|\int_{t_0}^{t_1}f(s,X_s)\dif s\Big\|_{L^p(\Om)}
    \lesssim_C \om_{**}(t_1-t_0)\|f\|_{\mL_T^\infty\bC^{-1,**}_{m,\eps}},\no
\end{align}
and
\begin{align}
    \Big\|\mE\Big(\int_{t_0}^{t_1}f(s,X_s)\dif s\,\big|\,\cF_{t_0}\Big)\Big\|_{L^p(\Om)}
    \lesssim_C \bar\om_{**}(t_1-t_0)\|f\|_{\mL_T^\infty\bC^{-1,**}_{m,\eps}},\no
\end{align}
where $\om_{**}(\de):=\sqrt{\de}/\ell_{m+2}^{\eps}(\de^{-1})$ and $\bar\om_{**}(\de):=\sqrt{\de}/\big((\prod_{i=1}^{m+1}\ell_i(\de^{-1}))\ell_{m+2}^{1+\eps}(\de^{-1})\big)$.

Moreover, for each $x\in\mR^d$, let $X^x$ denote the weak solution to \eqref{SDE}
with $X^x_0=x$. Then, for every $t\in(0,T]$, the law of $X^x_t$ admits a density
$\rho_t(x,\cdot)$ with respect to the Lebesgue measure and satisfies the following two-sided Gaussian estimate: there exist constants $C_0,C_1,\gamma_0,\gamma_1>0$, depending only on $T,d,m,\eps$ and $\|b\|_{\mL_T^\infty\bC^{-1,*}_{m,\eps}}$, such that for all $t\in(0,T]$, $x\in\mR^d$ and a.e. $y\in\mR^d$,
\begin{align}
    C_0t^{-\frac d2}\exp\Big(-\gamma_0\frac{|x-y|^2}{t}\Big)
    \leq \rho_t(x,y)
    \leq
    C_1t^{-\frac d2}\exp\Big(-\gamma_1\frac{|x-y|^2}{t}\Big).
    \no
\end{align}
\et

\br
The theorem imposes no moment condition and no density assumption on the initial distribution $\mu$. In particular, Dirac initial laws are included.
\er

\br
The regularity assumption $b\in\mL_T^\infty\bC^{-1,*}_{m,\eps}$ should be viewed as a logarithmically corrected version of the scaling-critical condition $b\in\mL_T^\infty\bC^{-1}$. The additional logarithmic strength is used only to obtain summability at the endpoint and the small factor in the Schauder estimate.
\er

\subsection{Background and motivation}

Stochastic differential equations with singular drifts have been studied from several complementary viewpoints. In the function-valued case, Veretennikov \cite{Ver79} proved strong existence and pathwise uniqueness for non-degenerate Brownian SDEs with bounded measurable drifts, showing that Brownian noise can restore uniqueness beyond the deterministic ODE regime. Krylov and R\"ockner \cite{KR05} treated time-dependent drifts $b\in\mL_T^q L^p$ under the subcritical condition $\ff d p+\ff2q<1$ and established strong well-posedness for the corresponding It\^o SDE. The Zvonkin--Krylov method was further developed for Sobolev diffusion coefficients and singular drifts by Zhang \cite{Zha05,Zha11,Zha16}, where strong solutions, stochastic homeomorphism flows and stability estimates were obtained. At the scaling-critical level $\ff d p+\ff2q=1$, Beck, Flandoli, Gubinelli and Maurelli \cite{BFGM19} studied stochastic ODEs and the associated linear PDE with critical drifts, proving regularity, duality and uniqueness results in critical classes; see also Krylov \cite{Kry21} for diffusion processes with drift in $L^d$ and R\"ockner and Zhao \cite{RZ23} for weak solutions with critical time-dependent drifts. For genuinely distribution-valued drifts, Bass and Chen \cite{BC01,BC03} initiated the theory by interpreting the drift through Dirichlet processes and additive functionals. Flandoli, Issoglio and Russo \cite{FIR17} introduced virtual solutions for multidimensional SDEs with time-independent drifts in negative H\"older spaces, typically $b\in\bC^\a$ with $\a>-\ff12$. Rough path and paracontrolled approaches provide another way to treat products involving distributions: Delarue and Diel \cite{DD16} handled one-dimensional time-dependent distributional drifts motivated by polymer models, Cannizzaro and Chouk \cite{CC18} constructed multidimensional solutions and polymer measures with white-noise potentials, and Kremp and Perkowski \cite{KP22} extended related ideas to SDEs driven by L\'evy noise. Zhang and Zhao \cite{ZZ18} developed a weak well-posedness theory for singular Brownian diffusions with $b\in H^\a_p$, $\a\in[-\ff12,0]$ and $p>\ff d{1+\a}$, where the drift integral is recovered through Krylov estimates and becomes a zero-energy additive functional. In the divergence-free direction, the structural cancellation has also been exploited through PDE methods: Seregin, Silvestre, \v{S}ver\'ak and Zlato\v{s} \cite{SSSZ12} studied regularity effects of divergence-free drifts, while Zhang and Zhao \cite{ZZ21} and Zhang \cite{Zha24} used maximum-principle and De Giorgi type arguments to construct weak solutions in supercritical regimes. More recently, energy methods have made it possible to go below the classical product threshold. Gr\"afner and Perkowski \cite{GP24} studied singular SDEs with divergence-free supercritical distributional drifts and obtained weak well-posedness of energy solutions under suitable density assumptions on the initial law. Hao and Zhang \cite{HZ25} considered divergence-free drifts $b\in\mL_T^q H^\a_p$ with $\a\in[-1,0]$, $p,q\in[2,\infty]$ and initial distributions with $L^2$ densities. Compared with these works, the present paper focuses on the endpoint critical class $\mL_T^\infty\bC^{-1}$ and shows that a logarithmic improvement, together with the divergence-free structure, is sufficient to recover weak well-posedness for arbitrary initial laws.

The scaling analysis explains why the present paper is concerned with a borderline problem. Formally, if $X$ solves \eqref{SDE}, then for $\lambda>0$ we set
\begin{align}
    X_t^\lambda:=\lambda^{-1}X_{\lambda^2t},\ 
    W_t^\lambda:=\lambda^{-1}W_{\lambda^2t},\ 
    b^\lambda(t,x):=\lambda b(\lambda^2t,\lambda x).\no
\end{align}
Then $X^\lambda$ formally satisfies
\begin{align}
    \dif X_t^\lambda=b^\lambda(t,X_t^\lambda)\dif t+\sqrt{2}\dif W_t^\lambda.\no
\end{align}
For the homogeneous H\"older--Besov scale one has
\begin{align}
    \|b^\lambda\|_{\mL^q(\mR_+;\dot{\bC}^{\a})}
    =\lambda^{1+\a-\frac2q}\|b\|_{\mL^q(\mR_+;\dot{\bC}^{\a})}.\no
\end{align}
Thus the regimes are classified as
\begin{align}
    \hbox{Subcritical: } \ff2q<1+\a;\ 
    \hbox{Critical: } \ff2q=1+\a;\ 
    \hbox{Supercritical: } \ff2q>1+\a.\no
\end{align}
In particular, the time-bounded drift class $\mL_T^\infty\bC^{-1}$ is exactly critical. The recent work \cite{HZ25} develops a systematic theory for first-order SDEs with distributional drifts in Bessel potential spaces. In the subcritical case, weak well-posedness is obtained by combining the Schauder estimates, Krylov estimates and Zvonkin transform. In the critical and supercritical divergence-free case, weak solutions and Markov selections can only be constructed by energy methods when the initial law has an $L^2$ density \cite{GP24,HZ25}. Our aim is to recover enough endpoint regularization at the critical order $-1$ by imposing a logarithmic correction, and thereby prove weak well-posedness for arbitrary initial laws.

The analytic obstruction is already visible from the Kolmogorov equation associated with \eqref{SDE}:
\begin{align}\label{intro-kolmogorov}
    \p_tu^\l=\Delta u^\l-\l u^\l+b\cdot\nn u^\l+f,
    \ u_0^\l\equiv0.
\end{align}
If $b\in\bC^\a$ with $\a<0$, the classical Schauder estimate gives at most $u^\l\in\bC^{2+\a}$. Hence $\nn u^\l$ has regularity $1+\a$, and the classical paraproduct criterion for the product $b\cdot\nn u^\l$ would require $\a+(1+\a)>0$, that is, $\a>-\frac12$. The endpoint $\a=-1$ is therefore far outside the classical product regime. This is the main reason why the usual Schauder--Zvonkin argument cannot be applied directly at the scaling-critical level.

The divergence-free assumption provides the structural cancellation needed to formulate the singular drift term. Using the Bony decomposition, for smooth $u$ one may write
\begin{align}
    b\cdot\nn u
    =\div(b\preceq u)+b\succ\nn u-\div b\preceq u
    =:b\odot\nn u-\div b\preceq u.\no
\end{align}
Thus, when $\div b=0$, the drift term is interpreted as $b\cdot\nn u=b\odot\nn u$; the precise definition is given in \eqref{decomposition}. The logarithmically corrected spaces $\bC^{-1,*}_{m,\eps}$ and $\bC^{1,**}_{m,\eps}$ are designed so that the product is meaningful in the following sense:
\begin{align}
    b\in\bC^{-1,*}_{m,\eps},\ u\in\bC^{1,**}_{m,\eps}
    \ \Rightarrow\ 
    b\cdot\nn u\in\bC^{-1,*}_{m,\eps}.\no
\end{align}
The logarithmic correction also supplies the small factor that replaces the missing positive power in the subcritical theory. More precisely, the heat resolvent gains two derivatives from the stronger $*$ scale to the weaker $**$ scale with an additional factor $\ff1{\ell_{m+2}(J_\l)}$, where $J_\l$ is defined in \eqref{Jl-def}. This logarithmic smallness allows us to solve \eqref{intro-kolmogorov} for large $\l$ and to choose the Zvonkin map $\Phi(t,x)=x+u^\l(t,x)$ with small gradient. On the probabilistic side, the same gain yields logarithmic Krylov estimates for distributional test functions, which makes it possible to construct the additive functional $A^b$ and to prove the stochastic substitution formula needed in the Zvonkin transformation. In this way, the present paper can be viewed as a logarithmically corrected critical counterpart of the subcritical theory in \cite{HZ25}, with the additional feature that no moment condition and no density assumption is imposed on the initial distribution $\mu$. Moreover, the divergence-free approximation also allows us to propagate uniform Aronson-type bounds to the singular limit, yielding the two-sided Gaussian estimate in Theorem \ref{main-result}.

The main contributions of this paper can be summarized as follows.
\begin{itemize}
\item We exploit a carefully chosen gap between two logarithmically weighted H\"older--Besov scales adapted to the critical drift regularity. This gap yields the logarithmic smallness needed to close the endpoint Schauder estimate.
\item We establish uniform logarithmic Krylov estimates for the smooth approximating SDEs. These estimates allow us to construct additive functionals associated with distributional test functions and to identify the drift term after passing to the limit.
\item We prove a logarithmic stochastic substitution formula by stochastic sewing and combine it with a Zvonkin transformation. This yields weak uniqueness in the logarithmic Krylov class for divergence-free distributional drifts of order $-1$.
\item We represent the approximating generators in divergence form with bounded measurable uniformly elliptic leading coefficients and bounded lower-order terms by using the divergence-free structure. This yields Aronson-type two-sided Gaussian bounds for the smooth approximating heat kernels, uniformly in the approximation parameter. After weak uniqueness has identified the limiting law, these bounds pass to the limit and give the two-sided Gaussian estimate for the transition density of the original singular SDE.
\end{itemize}

\subsection{Structure of the paper}

The rest of the paper is organized as follows.

In Section \ref{sec2}, we collect the analytic preliminaries. We first introduce the logarithmically corrected H\"older--Besov spaces $\bC^{\a,*}_{m,\eps}$ and $\bC^{\a,**}_{m,\eps}$, whose dyadic weights are defined through the iterated logarithms $\ell_i$. We then recall the Bony paraproduct decomposition and prove that, under the divergence-free condition $\div b=0$, the singular product $b\cdot\nn u$ is well defined when $b\in\bC^{-1,*}_{m,\eps}$ and $u\in\bC^{1,**}_{m,\eps}$.

In Section \ref{sec3}, we prove the Schauder estimates for the Kolmogorov equation \eqref{PDE}. The key point is the logarithmic gain from the stronger scale $\bC^{-1,*}_{m,\eps}$ to the weaker scale $\bC^{1,**}_{m,\eps}$, quantified by the factor $\ell_{m+2}(J_\l)^{-1}$. As a consequence, for all sufficiently large $\l$, we obtain well-posedness of the resolvent equation $u_t^\l=\sI_t^\l(b\cdot\nn u^\l+f)$ together with the a priori estimate $u^\l\in\mL_T^\infty\bC^{1,**}_{m,\eps}$.

In Section \ref{sec-log-krylov}, we convert the analytic estimates into probabilistic estimates for the smooth approximating SDEs with drifts $b_n=b*\phi_n$. We establish uniform logarithmic Krylov estimates which allow us to construct the additive functionals $(A_t^f)_n$ for distributional test functions $f$ and to prove the logarithmic stochastic substitution formula, which is the main probabilistic ingredient for the Zvonkin transformation.

In Section \ref{sec5}, we prove Theorem \ref{main-result}. We first construct weak solutions to \eqref{SDE} by tightness of the approximating solutions $X^n$, identify the limiting drift as the additive functional $A^b$, and show that the constructed solutions belong to the logarithmic Krylov class. We then solve the backward equation \eqref{zvonkin-pde}, define the Zvonkin map $\Phi(t,x)=x+u^\l(t,x)$, and prove uniqueness in law within this class. Finally, we derive the two-sided Gaussian bounds by applying the Aronson estimates uniformly to the smooth divergence-form approximations and passing to the limit.

We conclude the introduction with the following notations. Throughout this paper, we use $C$ with or without subscripts to denote constants, whose values may vary from line to line. Also, we use $:=$ to indicate a definition. We write $A\lesssim_C B$ and $A\asymp_C B$ or simply without subscripts to mean that for some constant $C\geq1$, $A\leq CB$ and $C^{-1}B\leq A\leq CB$. Moreover,
for the convenience of the readers, we list some frequently used notations:

\begin{itemize}
\item $\mN$: set of positive integers; $\mN_0:=\{0\}\cup\mN$.

\item $\mR^n$: Euclidean space of dimension $n\in\mN$.

\item $\bC^{\a,*}_{m,\eps},\ \bC^{\a,**}_{m,\eps}$: logarithmically corrected H\"older-Besov spaces (see Definition \ref{Holder} below).

\item $C_T:=C([0,T];\mR^d)$.

\item $P_t$: semigroup of the operator $\Delta$.

\item $f\prec g, f\circ g, f\succ g$: paraproduct. $f\preceq g:=f\prec g+f\circ g$ (see \eqref{bonyeq} below).

\end{itemize}

\section{Preliminaries}\label{sec2}
This section sets up the analytic framework used throughout the paper. We first introduce the logarithmically corrected H\"older-Besov spaces via a Littlewood-Paley decomposition. We then recall the paraproduct calculus needed to interpret the singular product $b\cdot\nabla u$ when $b$ has regularity $-1$.

Let $\sS(\mR^d)$ be the Schwartz space of all rapidly decreasing functions on $\mR^{d}$, and let $\sS'(\mR^{d})$ denote the dual space of $\sS(\mR^{d})$, the space of tempered distributions. For any $f\in\sS(\mR^{d})$, we define the Fourier transform $\hat f$ and inverse Fourier transform $\check f$ respectively by
\begin{align}
    \hat f(\xi):=\ff{1}{(2\pi)^{d/2}}\int_{\mR^{d}}\e^{-i\xi\cdot x}f(x)\dif x,\ \xi\in\mR^{d},\no
\end{align}
\begin{align}
    \check f(x):=\ff{1}{(2\pi)^{d/2}}\int_{\mR^{d}}\e^{i\xi\cdot x}f(\xi)\dif\xi, \ x\in\mR^{d}.\no
\end{align}
For $q\in[1,\infty]$ and a normed space $\mB$, we will use the notation
\begin{align}
    \mL_T^q\mB:=L^q([0,T];\mB),\ \mL_T^q:=\mL_T^qL^q,\no
\end{align}
where $(L^p,\|\cdot\|_p)$ is the common $L^p$ space in $\mR^d$. Additionally, we denote by $C_b=C_b(\mR^{d})$ (resp. $C_c=C_c(\mR^d)$) the Banach space of all bounded continuous functions on $\mR^{d}$ (resp. the space of all continuous functions on $\mR^{d}$ with compact support), and by $C_b^\infty=C_b^\infty(\mR^{d})$ (resp. $C_c^\infty=C_c^\infty(\mR^{d})$) the space of all smooth functions with bounded derivatives of all orders (resp. the space of all smooth functions with compact support).

\subsection{Logarithmically corrected H\"older-Besov spaces}\label{sec21}

In this subsection we recall the Littlewood--Paley definition of H\"older--Besov spaces and introduce their logarithmically corrected variants.
Let $\chi_{0}$ be a symmetric $C^\infty$-function on $\mR^{d}$ such that
$$\chi_{0}(\xi)=1 \text{\ for\ } \xi\in B_{1} \text{\ and\ }\chi_{0}(\xi)=0 \text{\ for\ } \xi\notin B_{4/3}.$$
For $j\in\mN$, define 
\begin{align}
    \phi_j(\xi):=\left\{
    \begin{aligned}
        &\chi_0(2^{-j}\xi)-\chi_0(2^{-(j-1)}\xi),\ &j\geq1,\\
        &\chi_0(\xi),&j=0.
    \end{aligned}
    \right.\no
\end{align}
This definition ensures that 
\begin{align}
    \sum_{j\geq0}\phi_j(\xi)=1,\ \forall\xi\in\mR^{d},\no
\end{align}
and
\begin{align}
    \supp(\phi_j)\subset\{\xi:2^{j-1}\leq|\xi|\leq2^{j+2}/3\},\ j\geq1,\ \supp(\phi_0)\subset B_{4/3}.\no
\end{align}
For given $j\geq0$, we define the dyadic block operator $\cR_j$ on $\sS'$ as
\begin{align}
    \cR_j f(x):=(\phi_j\hat f)\check\ (x)=\check\phi_j*f(x),\no
\end{align}
where the convolution is taken in the distributional sense. By scaling, we have
\begin{align}
    \check\phi_j(x)=2^{(j-1)d}\check\phi_1(2^{(j-1)}x),\ j\geq1.\no
\end{align}
For $j\in\mN$, it follows that
\begin{align}\label{wtR}
    \cR_j=\cR_j\wt\cR_j,\ \text{where}\ \wt\cR_j:=\sum_{|i-j|\leq1}\cR_i,
\end{align}
with the convention that $\cR_j:=0$ for $j<0$.

The following Bernstein inequality is standard.

\bl\label{bernstein}
For any $k\in\mN$ and $1\leq p\leq q\leq\infty$, there exists a constant $C=C(k,p,q,d)>0$ such that for all $j\geq0$, 
\begin{align}
    \|\nn^{k}\cR_j f\|_{q}\lesssim_C 2^{j(k+d(\ff{1}{p}-\ff{1}{q}))}\|\cR_j f\|_{p}.\no
\end{align}
\el

Now we introduce the logarithmically corrected H\"older-Besov spaces used below. For $r\geq0$, define
\begin{align}\label{ellm}
\ell_0(r):=1+r,\ \ell_1(r):=\log(e+r),\ \ell_{k+1}(r):=\log(e+\ell_k(r)),\ k\ge1.
\end{align}
For fixed $\eps>0$, $m\in\mN$ and $\theta\in\{0,1\}$, define
$$w_{m,\eps,\theta}(j):=\Big(\prod_{i=0}^m \ell_i(j)\Big)\ell_{m+1}^{1+\eps}(j)\ell_{m+2}^{\theta}(j),\ j\geq0.$$
Compared with the logarithmic Besov framework developed in \cite{CDFM15}, we need a slight refinement. The parameter $\theta$ distinguishes the two dyadic weights used below: the case $\theta=1$ contains one additional iterated logarithmic gain compared with the case $\theta=0$.

\bd\label{Holder}
For any $\a\in\mR$, the usual H\"older--Besov space is defined by
$$\bC^\a:=\Big\{f\in\sS':\|f\|_{\bC^\a}:=\sup_{j\geq0}2^{j\a}\|\cR_j f\|_{\infty}<\infty\Big\}.$$
Let $\eps>0$ and $m\in\mN$ be fixed. For $\theta\in\{0,1\}$, we define the logarithmically corrected H\"older space
$$\bC^{\a}_{m,\eps,\theta}:=\Big\{f\in\sS':\|f\|_{\bC^{\a}_{m,\eps,\theta}}:=\sup_{j\geq0}w_{m,\eps,\theta}(j)2^{j\a}\|\cR_j f\|_{\infty}<\infty\Big\}.$$
For later use, we denote
$$\bC^{\a,*}_{m,\eps}:=\bC^{\a}_{m,\eps,1},\ \bC^{\a,**}_{m,\eps}:=\bC^{\a}_{m,\eps,0}.$$
\ed

\br
The space $\bC^{\a,*}_{m,\eps}$ is slightly stronger than $\bC^{\a,**}_{m,\eps}$. On the frequency side $|\xi|\sim2^j$, the factor $\ell_0(j)=1+j$ corresponds to the usual logarithmic correction $\log(e+|\xi|)$, while $\ell_k(j)$, $k\geq1$, corresponds to higher iterated logarithmic corrections. Thus these spaces should be viewed as H\"older-Besov spaces with iterated logarithmic gains.
\er
\br
For $\a=0$, $\bC^{0,*}_{m,\eps}\subset \bC^{0,**}_{m,\eps}\subset C_b.$ We emphasize that the condition $\eps>0$ is essential here. Indeed, the strict exponent $1+\eps$ is precisely what makes the logarithmic tail summable; without this, it would in general fail to imply the embedding into $C_b$.
\er

\subsection{Paraproduct calculus}
To ensure that the term $b\cdot\nn u$ makes sense when $b$ is distribution-valued and $u\in\sS(\mR^{d})$, we recall the Bony paraproduct decomposition, which allows us to control the product of functions in H\"{o}lder spaces under weaker regularity assumptions. For any $k\geq0$, define the cut-off low frequency operator $S_k$ as
\begin{align}
    S_k f:=\sum_{j=0}^{k-1}\cR_j f\overset{k\rightarrow\infty}{\rightarrow}f,\no
\end{align} 
where we take $S_{-1}=0$ by convention. For $f,g\in\sS'(\mR^{d})$, we define the following paraproducts
\begin{align}\label{bonyeq}
    f\prec g:=\sum_{k\geq0}S_{k-1} f\cR_k g,\ f\circ g:=\sum_{k\geq0}\cR_k f\wt\cR_k g,
\end{align}
where $\wt\cR_j$ is defined in \eqref{wtR}. The Bony decomposition of $fg$ is formally written as (cf. \cite{BCD11})
\begin{align}
    fg=f\prec g+f\circ g+g\prec f=:f\preceq g+f\succ g.\no
\end{align} 
The essential property of the Bony decomposition is that
\begin{align}
    \cR_k(S_{j-1}f\cR_j g)=0\ \text{for}\ |k-j|>2.\no
\end{align}

Using the Bony decomposition we can write
\begin{align}\label{decomposition}
    b\cdot\nabla u:={\div(b\preceq u)+b\succ\nabla u}-\div b\preceq u=:b\odot\nabla u-\div b\preceq u.
\end{align}
In particular, if $\div b=0$, then 
\begin{align}
    b\cdot\nabla u=b\odot\nabla u.\no
\end{align}

We have the following result.

\bl\label{driftterm}
Assume that $b\in \bC^{-1,*}_{m,\eps}$ is divergence free and $u\in \bC^{1,**}_{m,\eps}$. Then $b\cdot\nn u\in \bC^{-1,*}_{m,\eps}$.
\el

\bpf
For fixed $j\geq0$, we have
\begin{align}
   \|\cR_j(b\prec u)\|_{\infty}&=\Big\|\cR_j\Big(\sum_{k\geq0}S_{k-1} b\cR_k u\Big)\Big\|_{\infty}\no\\
   &\lesssim \sum_{|k-j|\leq2}\|S_{k-1} b\|_{\infty}\|\cR_k u\|_{\infty}\lesssim\|S_{j-1}b\|_{\infty}\|\cR_j u\|_{\infty}\no\\
   &\lesssim\Big(\sum_{\ell=0}^{j-2}\frac{2^{\ell}\|b\|_{\bC^{-1,*}_{m,\eps}}}{w_{m,\eps,1}(\ell)}\Big)\frac{2^{-j}\|u\|_{\bC^{1,**}_{m,\eps}}}{w_{m,\eps,0}(j)}\lesssim\frac{\|b\|_{\bC^{-1,*}_{m,\eps}}\|u\|_{\bC^{1,**}_{m,\eps}}}{w_{m,\eps,1}(j)}.\no
\end{align}
Similarly, for fixed $j\geq0$,
\begin{align}
   \|\cR_j(b\circ u)\|_{\infty}&=\Big\|\cR_j\Big(\sum_{k\geq0}\cR_k b\wt\cR_k u\Big)\Big\|_{\infty}\lesssim \sum_{k\geq j-3}\|\cR_k b\|_{\infty}\|\cR_k u\|_{\infty}\no\\
   &\lesssim\sum_{k\geq j}\Big(\frac{\|b\|_{\bC^{-1,*}_{m,\eps}}}{w_{m,\eps,1}(k)}\frac{\|u\|_{\bC^{1,**}_{m,\eps}}}{w_{m,\eps,0}(k)}\Big)\lesssim\frac{\|b\|_{\bC^{-1,*}_{m,\eps}}\|u\|_{\bC^{1,**}_{m,\eps}}}{w_{m,\eps,1}(j)}.\no
\end{align}
Hence, by the Bernstein inequality, we have $\div(b\preceq u)\in \bC^{-1,*}_{m,\eps}$. For $b\succ\nn u$, we have
\begin{align}
   \|\cR_j(b\succ \nn u)\|_{\infty}&=\Big\|\cR_j\Big(\sum_{k\geq0}\cR_k b S_{k-1}(\nn u)\Big)\Big\|_{\infty}\no\\
   &\lesssim \sum_{|k-j|\leq2}\|\cR_k b\|_{\infty}\|S_{k-1}(\nn u)\|_{\infty}\lesssim\|\cR_j b\|_{\infty}\|S_{j-1}(\nn u)\|_{\infty}\no\\
   &\lesssim\frac{2^j\|b\|_{\bC^{-1,*}_{m,\eps}}}{w_{m,\eps,1}(j)}\Big(\sum_{\ell=0}^{\infty}\frac{\|u\|_{\bC^{1,**}_{m,\eps}}}{w_{m,\eps,0}(\ell)}\Big)\lesssim\frac{2^j\|b\|_{\bC^{-1,*}_{m,\eps}}\|u\|_{\bC^{1,**}_{m,\eps}}}{w_{m,\eps,1}(j)},\no
\end{align}
which shows $b\succ\nn u\in \bC^{-1,*}_{m,\eps}$. By \eqref{decomposition}, the proof is complete.
\epf

\br\label{bvariant}
A similar calculation shows that for $b\in \bC^{-1}$ divergence free and $u\in\bC^{1,**}_{m,\eps}$, we have $b\cdot\nn u\in \bC^{-1}$. This explains the particular choice of the logarithmic correction in the definition of the space $\bC^{1,**}_{m,\eps}$.
\er

\section{Schauder estimate}\label{sec3}
The paraproduct estimate from the previous section allows us to define the drift term in the Kolmogorov equation. We now prove the Schauder estimates that will later provide the logarithmic small factor needed for the Zvonkin transform. In this section we consider the equation
\begin{align}\label{PDE}
    \p_t u^\l=\Delta u^\l-\l u^\l+b\cdot\nn u^\l+f,\ u^\l_0\equiv0.
\end{align}
Let $(P_t)_{t\geq0}$ be the Gaussian heat semigroup given by
\begin{align}
    P_t f=\varphi_t*f,\ \varphi_t(x)=(4\pi t)^{-\frac{d}{2}}\e^{-\frac{|x|^2}{4t}}.\no
\end{align}
For $\l\geq0$ and a time-dependent distribution $f_t(\cdot):\mR_+\to\sS'(\mR^{d})$, we define
\begin{align}
    u_t^\l(z):=\sI_t^\l(f)(z):=\int_0^t \e^{-\l(t-s)}P_{t-s}f_s(z)\dif s,\ t>0.\no
\end{align}
Then, in the sense of distributions,
\begin{align}
    \partial_t u^\l=(\Delta-\l)u^\l+f.\no
\end{align}
For $\l\geq1$, set
\begin{align}\label{Jl-def}
    J_\l:=2+\Big\lfloor \ff12\log_2(\e+\l)\Big\rfloor .
\end{align}
Then, for each fixed $i\geq1$,
\begin{align}\label{Jl-asymp}
    J_\l\asymp \log(\e+\l),\ \ell_i(J_\l)\asymp \ell_{i+1}(\l).
\end{align}

We first record the dyadic estimate which will be used in the Schauder estimate below.

\bl
Let $T>0$ and $\l\geq1$. For every $j\geq0$,
\begin{align}\label{block-resolvent}
    \|\cR_j\sI^\l(f)\|_\infty
    \lesssim 2^{-2j}\Big(1\wedge\ff{2^{2j}}{\l}\Big)
    \|\cR_j f\|_\infty.
\end{align}
\el

\bpf
By \eqref{wtR},
\begin{align}
    \|\cR_j\sI_t^\l(f)\|_\infty
    &\leq\int_0^t \e^{-\l(t-s)}\|\cR_j(\varphi_{t-s}*f_s)\|_\infty\dif s\no\\
    &=\int_0^t \e^{-\l(t-s)}\|\wt\cR_j\varphi_{t-s}*\cR_jf_s\|_\infty\dif s\no\\
    &\lesssim \|\cR_j f\|_\infty
    \int_0^T \e^{-\l s}\|\wt\cR_j\varphi_s\|_1\dif s.\no
\end{align}
By \cite{HWZ20}, we have $\|\wt\cR_j\varphi_s\|_1\lesssim 1\wedge(2^{2j}s)^{-2}$. Consequently,
\begin{align}
    \int_0^T \e^{-\l s}\|\wt\cR_j\varphi_s\|_1\dif s&\lesssim \int_0^\infty \e^{-\l s}\big[1\wedge(2^{2j}s)^{-2}\big]\dif s\no\\
    &=2^{-2j}\int_0^\infty \e^{-\l 2^{-2j} r}\big[1\wedge r^{-2}\big]\dif r\lesssim 2^{-2j}\Big(1\wedge\ff{2^{2j}}{\l}\Big),\no
\end{align}
which gives \eqref{block-resolvent}.
\epf

We shall also use the following elementary consequences.

\bl\label{Jl-sum-es}
For fixed $\eps>0$, $m\in\mN$ and $\l\geq1$,
\begin{align}\label{sum-schauder}
    \sup_{j\geq0}\ff{w_{m,\eps,0}(j)}{w_{m,\eps,1}(j)}\Big(1\wedge\ff{2^{2j}}{\l}\Big)
    \lesssim \ff1{\ell_{m+2}(J_\l)},
\end{align}
\begin{align}\label{sum-star}
    \sum_{j\geq0}\ff{2^{-j}}{w_{m,\eps,1}(j)}\Big(1\wedge\ff{2^{2j}}{\l}\Big)
    \lesssim \ff{\l^{-1/2}}{w_{m,\eps,1}(J_\l)},\ 
    \sum_{j\geq0}\ff1{w_{m,\eps,1}(j)}\Big(1\wedge\ff{2^{2j}}{\l}\Big)
    \lesssim \ff1{\ell_{m+1}^{\eps}(J_\l)\ell_{m+2}(J_\l)},
\end{align}
\begin{align}\label{sum-starstar}
    \sum_{j\geq0}\ff{2^{-j}}{w_{m,\eps,0}(j)}\Big(1\wedge\ff{2^{2j}}{\l}\Big)
    \lesssim \ff{\l^{-1/2}}{w_{m,\eps,0}(J_\l)},\ 
    \sum_{j\geq0}\ff1{w_{m,\eps,0}(j)}\Big(1\wedge\ff{2^{2j}}{\l}\Big)
    \lesssim \ff1{\ell_{m+1}^{\eps}(J_\l)}.
\end{align}
\el

\bpf
Since $2^{J_\l}\asymp \l^{1/2}$, the proof is obtained by splitting the sums at $J_\l$. For $\theta\in\{0,1\}$, by the monotonicity and slow variation of $w_{m,\eps,\theta}$,
\begin{align}
    \sum_{j\geq2}\ff{2^{-j}}{w_{m,\eps,\theta}(j)}\Big(1\wedge\ff{2^{2j}}{\l}\Big)\lesssim 
    \ff1\l\sum_{2\leq j\leq J_\l}\ff{2^j}{w_{m,\eps,\theta}(j)}+\sum_{j>J_\l}\ff{2^{-j}}{w_{m,\eps,\theta}(j)}
    \lesssim \ff{\l^{-1/2}}{w_{m,\eps,\theta}(J_\l)}.\no
\end{align}
This gives the first estimates in \eqref{sum-star} and \eqref{sum-starstar}. Similarly,
\begin{align}
    \sum_{j\geq2}\ff1{w_{m,\eps,\theta}(j)}\Big(1\wedge\ff{2^{2j}}{\l}\Big)
    \lesssim \ff1{w_{m,\eps,\theta}(J_\l)}+
    \sum_{j>J_\l}\ff1{w_{m,\eps,\theta}(j)}.\no
\end{align}
By the integral test,
\begin{align}
    \sum_{j>J_\l}\ff1{w_{m,\eps,0}(j)}\lesssim
    \ff1{\ell_{m+1}^{\eps}(J_\l)},\ 
    \sum_{j>J_\l}\ff1{w_{m,\eps,1}(j)}\lesssim
    \ff1{\ell_{m+1}^{\eps}(J_\l)\ell_{m+2}(J_\l)}.\no
\end{align}
Hence the second estimates in \eqref{sum-star} and \eqref{sum-starstar} follow. Finally, since
$w_{m,\eps,0}(j)/w_{m,\eps,1}(j)=1/\ell_{m+2}(j)$, splitting into $j\leq J_\l$ and $j>J_\l$ gives
\begin{align}
    \sup_{j\geq2}\ff{w_{m,\eps,0}(j)}{w_{m,\eps,1}(j)}\Big(1\wedge\ff{2^{2j}}{\l}\Big)
    \lesssim \ff1{\ell_{m+2}(J_\l)},\no
\end{align}
which proves \eqref{sum-schauder}.
\epf

\br
The bounds in Lemma \ref{Jl-sum-es} are actually equivalence estimates. Namely, each $\lesssim$ there can be replaced by $\asymp$, with constants independent of $\lambda\geq 1$ but depending on $m$ and $\eps$.
\er

To obtain the well-posedness of \eqref{PDE}, we need the following Schauder estimate.

\bl\label{schauder-estimate}
Let $\a\in\mR$, $\eps>0$ and $m\in\mN$. For any $T>0$, there is a constant $C=C(T,d,m,\eps)$ such that for every $\l\geq1$,
\begin{align}\label{schauder1}
    \|\sI_t^\l(f)\|_{\mL_T^\infty \bC^{\a,**}_{m,\eps}}\lesssim_C \|f\|_{\mL_T^\infty\bC^{\a-2,**}_{m,\eps}},
\end{align}
and
\begin{align}\label{schauder2}
    \|\sI_t^\l(f)\|_{\mL_T^\infty \bC^{\a,**}_{m,\eps}}
    \lesssim_C \ff1{\ell_{m+2}(J_\l)}\|f\|_{\mL_T^\infty\bC^{\a-2,*}_{m,\eps}}.
\end{align}
\el

\bpf
For \eqref{schauder1}, by \eqref{block-resolvent}, for every $j\geq0$,
\begin{align}
    w_{m,\eps,0}(j)2^{j\a}\|\cR_j\sI^\l(f)\|_{\mL_T^\infty L^\infty}
    \lesssim \Big(1\wedge\ff{2^{2j}}{\l}\Big)
    \|f\|_{\mL_T^\infty\bC^{\a-2,**}_{m,\eps}}
    \lesssim \|f\|_{\mL_T^\infty\bC^{\a-2,**}_{m,\eps}}.\no
\end{align}
Taking the supremum over $j$ gives \eqref{schauder1}.

It remains to prove \eqref{schauder2}. By \eqref{block-resolvent},
\begin{align}
    w_{m,\eps,0}(j)2^{j\a}\|\cR_j\sI^\l(f)\|_{\mL_T^\infty L^\infty}
    &\lesssim \ff{w_{m,\eps,0}(j)}{w_{m,\eps,1}(j)}
    \Big(1\wedge\ff{2^{2j}}{\l}\Big)
    \|f\|_{\mL_T^\infty\bC^{\a-2,*}_{m,\eps}}.\no
\end{align}
Taking the supremum over $j$ and using \eqref{sum-schauder}, we obtain \eqref{schauder2}.
\epf

We have the following main result of this section.

\bt
Assume that $\eps>0$, $m\in\mN$, $\|b\|_{\mL_T^\infty\bC^{-1,*}_{m,\eps}}<\infty$ and $\div b=0$. Then there exists a constant
$\l_0=\l_0(T,d,m,\eps,\|b\|_{\mL_T^\infty\bC^{-1,*}_{m,\eps}})>0$, such that for any $\l\geq\l_0$ and
$f\in\mL_T^\infty\bC^{-1,**}_{m,\eps}$, there exists a unique $u^\l$ solving the equation
\begin{align}
    u_t^\l=\sI_t^\l(b\cdot\nn u^\l+f),\no
\end{align}
where the term $b\cdot\nn u^\l$ is defined in the sense of Lemma \ref{driftterm}, and such that
\begin{align}\label{pdees}
    \|u^\l\|_{\mL_T^\infty\bC^{1,**}_{m,\eps}}\lesssim_C\|f\|_{\mL_T^\infty\bC^{-1,**}_{m,\eps}}.
\end{align}
\et

\bpf
We only prove the a priori estimate \eqref{pdees}. The existence follows a standard Picard iteration. By Lemma \ref{driftterm}, Schauder estimate \eqref{schauder1} and \eqref{schauder2}, for suitably large $\l$,
\begin{align}
    \|u^\l\|_{\mL_T^\infty\bC^{1,**}_{m,\eps}}&=\|\sI_t^\l(b\cdot\nn u^\l+f)\|_{\mL_T^\infty\bC^{1,**}_{m,\eps}}\no\\
    &\lesssim \ff1{\ell_{m+2}(J_\l)}
    \|b\cdot\nn u^\l\|_{\mL_T^\infty\bC^{-1,*}_{m,\eps}}+
    \|f\|_{\mL_T^\infty\bC^{-1,**}_{m,\eps}}\no\\
    &\lesssim \ff1{\ell_{m+2}(J_\l)}
    \|b\|_{\mL_T^\infty\bC^{-1,*}_{m,\eps}}
    \|u^\l\|_{\mL_T^\infty\bC^{1,**}_{m,\eps}}+
    \|f\|_{\mL_T^\infty\bC^{-1,**}_{m,\eps}}.\no
\end{align}
Noting $\ell_{m+2}(J_\l)\to\infty$ as $\l\to\infty$, we complete the proof.
\epf

\br
Note that all terms on the right side of the equality are well-defined and the solution $u^\l$ constructed above satisfies that for any test function $\varphi\in C_c^\infty(\mR^{d})$, 
$$\p_t\<u^\l,\varphi\>=\<u^\l,\Delta\varphi\>-\l\<u^\l,\varphi\>+\<b\odot\nabla u^\l,\varphi\>+\<f,\varphi\>,$$
where $\<u^\l,\varphi\>:=\int_{\mR^{d}}u^\l(z)\varphi(z)\dif z$. This formulation indicates that $u$ is also a weak solution of PDE \eqref{PDE}.
\er

\section{Logarithmic Krylov estimate}\label{sec-log-krylov}

We next convert the Schauder estimates into probabilistic estimates along the smooth approximating solutions. These logarithmic Krylov estimates are the main tool for constructing additive functionals associated with distributional test functions. Throughout this section, $J_\l$ is defined by \eqref{Jl-def}, and we use the asymptotic relation \eqref{Jl-asymp}. We first prove the following estimates for the free resolvent.

\bl\label{endpoint-free}
Let $T>0$, $\eps>0$, $m\in\mN$ and $\l\geq1$. If $f\in\mL_T^\infty\bC^{-1,*}_{m,\eps}$, then
\begin{align}\label{free-star-Linfty}
    \|\sI^\l(f)\|_\infty\lesssim \ff{\l^{-1/2}}{w_{m,\eps,1}(J_\l)}\|f\|_{\mL_T^\infty\bC^{-1,*}_{m,\eps}},
\end{align}
\begin{align}\label{free-star-grad}
    \|\nn\sI^\l(f)\|_\infty\lesssim \ff1{\ell_{m+1}^{\eps}(J_\l)\ell_{m+2}(J_\l)}\|f\|_{\mL_T^\infty\bC^{-1,*}_{m,\eps}}.
\end{align}
Moreover, if $f\in\mL_T^\infty\bC^{-1,**}_{m,\eps}$, then
\begin{align}\label{free-starstar-Linfty}
    \|\sI^\l(f)\|_\infty\lesssim \ff{\l^{-1/2}}{w_{m,\eps,0}(J_\l)}\|f\|_{\mL_T^\infty\bC^{-1,**}_{m,\eps}},
\end{align}
\begin{align}\label{free-starstar-grad}
    \|\nn\sI^\l(f)\|_\infty\lesssim \ff1{\ell_{m+1}^{\eps}(J_\l)}\|f\|_{\mL_T^\infty\bC^{-1,**}_{m,\eps}}.
\end{align}
\el

\bpf
First assume $f\in\mL_T^\infty\bC^{-1,*}_{m,\eps}$. Note
\begin{align}
    \|\cR_j f\|_\infty
    \lesssim \ff{2^j}{w_{m,\eps,1}(j)}\|f\|_{\mL_T^\infty\bC^{-1,*}_{m,\eps}}.\no
\end{align}
Hence, by \eqref{block-resolvent} and \eqref{sum-star},
\begin{align}
    \|\sI^\l(f)\|_\infty\lesssim
    \sum_{j\geq0}\ff{2^{-j}}{w_{m,\eps,1}(j)}\Big(1\wedge\ff{2^{2j}}{\l}\Big)
    \|f\|_{\mL_T^\infty\bC^{-1,*}_{m,\eps}}
    \lesssim \ff{\l^{-1/2}}{w_{m,\eps,1}(J_\l)}\|f\|_{\mL_T^\infty\bC^{-1,*}_{m,\eps}},\no
\end{align}
and, by the Bernstein inequality,
\begin{align}
    \|\nn\sI^\l(f)\|_\infty
    \lesssim\sum_{j\geq0}\ff1{w_{m,\eps,1}(j)}\Big(1\wedge\ff{2^{2j}}{\l}\Big)
    \|f\|_{\mL_T^\infty\bC^{-1,*}_{m,\eps}}
    \lesssim \ff1{\ell_{m+1}^{\eps}(J_\l)\ell_{m+2}(J_\l)}
    \|f\|_{\mL_T^\infty\bC^{-1,*}_{m,\eps}}.\no
\end{align}
This proves \eqref{free-star-Linfty} and \eqref{free-star-grad}.

Next assume $f\in\mL_T^\infty\bC^{-1,**}_{m,\eps}$. Similarly, note
\begin{align}
    \|\cR_j f\|_\infty
    \lesssim \ff{2^j}{w_{m,\eps,0}(j)}\|f\|_{\mL_T^\infty\bC^{-1,**}_{m,\eps}}.\no
\end{align}
Therefore, by \eqref{block-resolvent} and \eqref{sum-starstar},
\begin{align}
    \|\sI^\l(f)\|_\infty\lesssim
    \sum_{j\geq0}\ff{2^{-j}}{w_{m,\eps,0}(j)}\Big(1\wedge\ff{2^{2j}}{\l}\Big)
    \|f\|_{\mL_T^\infty\bC^{-1,**}_{m,\eps}}
    \lesssim \ff{\l^{-1/2}}{w_{m,\eps,0}(J_\l)}\|f\|_{\mL_T^\infty\bC^{-1,**}_{m,\eps}},\no
\end{align}
\begin{align}
    \|\nn\sI^\l(f)\|_\infty\lesssim
    \sum_{j\geq0}\ff1{w_{m,\eps,0}(j)}\Big(1\wedge\ff{2^{2j}}{\l}\Big)
    \|f\|_{\mL_T^\infty\bC^{-1,**}_{m,\eps}}
    \lesssim \ff1{\ell_{m+1}^{\eps}(J_\l)}
    \|f\|_{\mL_T^\infty\bC^{-1,**}_{m,\eps}}.\no
\end{align}
This proves \eqref{free-starstar-Linfty} and \eqref{free-starstar-grad}. The proof is complete.
\epf

\bc\label{endpoint-kolmogorov}
Assume that $\eps>0$, $m\in\mN$, $\|b\|_{\mL_T^\infty\bC^{-1,*}_{m,\eps}}<\infty$ and $\div b=0$. There exists
$\l_0=\l_0(T,d,m,\eps,\|b\|_{\mL_T^\infty\bC^{-1,*}_{m,\eps}})>0$ such that for all $\l\geq\l_0$, the solution $u^\l$ to
\begin{align}
    u_t^\l=\sI_t^\l(b\cdot\nn u^\l+f)\no
\end{align}
satisfies the following estimates: If $f\in\mL_T^\infty\bC^{-1,*}_{m,\eps}$, then
\begin{align}\label{kol-star-Linfty}
    \|u^\l\|_\infty
    \lesssim \ff{\l^{-1/2}}{w_{m,\eps,1}(J_\l)}\|f\|_{\mL_T^\infty\bC^{-1,*}_{m,\eps}},
\end{align}
\begin{align}\label{kol-star-grad}
    \|\nn u^\l\|_\infty
    \lesssim \ff1{\ell_{m+1}^{\eps}(J_\l)\ell_{m+2}(J_\l)}\|f\|_{\mL_T^\infty\bC^{-1,*}_{m,\eps}}.
\end{align}
If $f\in\mL_T^\infty\bC^{-1,**}_{m,\eps}$, then
\begin{align}\label{kol-starstar-Linfty}
    \|u^\l\|_\infty
    \lesssim \ff{\l^{-1/2}}{w_{m,\eps,0}(J_\l)}\|f\|_{\mL_T^\infty\bC^{-1,**}_{m,\eps}},
\end{align}
\begin{align}
\label{kol-starstar-grad}
    \|\nn u^\l\|_\infty
    \lesssim \ff1{\ell_{m+1}^{\eps}(J_\l)}\|f\|_{\mL_T^\infty\bC^{-1,**}_{m,\eps}}.
\end{align}
\ec

\bpf
We first consider $f\in\mL_T^\infty\bC^{-1,*}_{m,\eps}$. By Lemma \ref{driftterm} and \eqref{schauder2},
\begin{align}
    \|u^\l\|_{\mL_T^\infty\bC^{1,**}_{m,\eps}}
    &\lesssim \ff1{\ell_{m+2}(J_\l)}
    \|b\cdot\nn u^\l+f\|_{\mL_T^\infty\bC^{-1,*}_{m,\eps}}\no\\
    &\lesssim \ff1{\ell_{m+2}(J_\l)}
    \|b\|_{\mL_T^\infty\bC^{-1,*}_{m,\eps}}
    \|u^\l\|_{\mL_T^\infty\bC^{1,**}_{m,\eps}}
    +\ff1{\ell_{m+2}(J_\l)}
    \|f\|_{\mL_T^\infty\bC^{-1,*}_{m,\eps}}.\no
\end{align}
Since $\ell_{m+2}(J_\l)\to\infty$, by taking $\l_0$ large enough we get
\begin{align}
    \|u^\l\|_{\mL_T^\infty\bC^{1,**}_{m,\eps}}
    \lesssim \ff1{\ell_{m+2}(J_\l)}
    \|f\|_{\mL_T^\infty\bC^{-1,*}_{m,\eps}}.\no
\end{align}
In particular,
\begin{align}
    \|b\cdot\nn u^\l\|_{\mL_T^\infty\bC^{-1,*}_{m,\eps}}
    \lesssim \ff1{\ell_{m+2}(J_\l)}
    \|b\|_{\mL_T^\infty\bC^{-1,*}_{m,\eps}}
    \|f\|_{\mL_T^\infty\bC^{-1,*}_{m,\eps}}.\no
\end{align}
Note $u^\l=\sI^\l(f)+\sI^\l(b\cdot\nn u^\l)$. By Lemma \ref{endpoint-free}, we obtain \eqref{kol-star-Linfty} and \eqref{kol-star-grad}.

Now let $f\in\mL_T^\infty\bC^{-1,**}_{m,\eps}$. By \eqref{pdees},
\begin{align}
    \|u^\l\|_{\mL_T^\infty\bC^{1,**}_{m,\eps}}
    \lesssim \|f\|_{\mL_T^\infty\bC^{-1,**}_{m,\eps}}.\no
\end{align}
Thus, by Lemma \ref{driftterm},
\begin{align}
    \|b\cdot\nn u^\l\|_{\mL_T^\infty\bC^{-1,*}_{m,\eps}}
    \lesssim \|b\|_{\mL_T^\infty\bC^{-1,*}_{m,\eps}}
    \|f\|_{\mL_T^\infty\bC^{-1,**}_{m,\eps}}.\no
\end{align}
By similar calculation and noting that
\begin{align}
    \ff{\l^{-1/2}}{w_{m,\eps,1}(J_\l)}
    \lesssim \ff{\l^{-1/2}}{w_{m,\eps,0}(J_\l)},\ 
    \ff1{\ell_{m+1}^{\eps}(J_\l)\ell_{m+2}(J_\l)}
    \lesssim\ff1{\ell_{m+1}^{\eps}(J_\l)},\no
\end{align}
we obtain \eqref{kol-starstar-Linfty} and \eqref{kol-starstar-grad}.
\epf

Let $b_n:=b*\phi_n\in\mL_T^\infty C_b^\infty$ be the smooth approximation of $b$. The approximation satisfies the uniform bound
\begin{align}
    \sup_n\|b_n\|_{\mL_T^\infty\bC^{-1,*}_{m,\eps}}
    \lesssim\|b\|_{\mL_T^\infty\bC^{-1,*}_{m,\eps}}<\infty.\no
\end{align}
Let $X^n$ solve the approximation equation:
\begin{align}
    \dif X_t^n=b_n(t,X_t^n)\dif t+\sqrt2\dif W_t,\ X_0^n\sim\mu\in\sP(\mR^d).\no
\end{align}
For $0<\de\leq T$, set
\begin{align}
    \om_*(\de):=\ff{\sqrt\de}{\ell_{m+2}^{\eps}(\de^{-1})\ell_{m+3}(\de^{-1})},\ 
    \om_{**}(\de):=\ff{\sqrt\de}{\ell_{m+2}^{\eps}(\de^{-1})},\no
\end{align}
and
\begin{align}
    \bar\om_*(\de):=\ff{\sqrt\de}
    {\big(\prod_{i=1}^{m+1}\ell_i(\de^{-1})\big)\ell_{m+2}^{1+\eps}(\de^{-1})\ell_{m+3}(\de^{-1})},\ 
    \bar\om_{**}(\de):=\ff{\sqrt\de}
    {\big(\prod_{i=1}^{m+1}\ell_i(\de^{-1})\big)\ell_{m+2}^{1+\eps}(\de^{-1})}.\no
\end{align}
We prove the following logarithmic Krylov estimate for $X^n$.

\bp\label{log-krylov}
Let $p\geq2$, $0\leq\tau_0\leq\tau_1\leq T$ be two stopping times with $\tau_1-\tau_0\leq\de$.

If $f\in\mL_T^1 C_b\cap\mL_T^\infty\bC^{-1,*}_{m,\eps}$, then there is a constant
$C_1:=C_1(p,T,d,m,\eps,\|b\|_{\mL_T^\infty\bC^{-1,*}_{m,\eps}})$ such that
\begin{align}\label{krylov-star}
    \sup_n\Big\|\int_{\tau_0}^{\tau_1}f(s,X_s^n)\dif s\Big\|_{L^p(\Om)}
    \lesssim_{C_1}\om_*(\de)\|f\|_{\mL_T^\infty\bC^{-1,*}_{m,\eps}}.
\end{align}
If $f\in\mL_T^1 C_b\cap\mL_T^\infty\bC^{-1,**}_{m,\eps}$, then there is a constant
$C_2:=C_2(p,T,d,m,\eps,\|b\|_{\mL_T^\infty\bC^{-1,*}_{m,\eps}})$ such that
\begin{align}\label{krylov-starstar}
    \sup_n\Big\|\int_{\tau_0}^{\tau_1}f(s,X_s^n)\dif s\Big\|_{L^p(\Om)}
    \lesssim_{C_2}
    \om_{**}(\de)\|f\|_{\mL_T^\infty\bC^{-1,**}_{m,\eps}}.
\end{align}
\ep

\bpf
First assume $f\in C_c^\infty([0,T]\times\mR^d)$. By time reversal, Corollary \ref{endpoint-kolmogorov} also applies to the backward equation
\begin{align}
    \p_t u_n^\l+\Delta u_n^\l+b_n\cdot\nn u_n^\l-\l u_n^\l+f=0,
    \ u_{n}^\l(T)=0.\no
\end{align}
By the It\^o formula,
\begin{align}
    \int_0^t f(s,X_s^n)\dif s
    =u_n^\l(0,X_0^n)-u_n^\l(t,X_t^n)
    +\l\int_0^tu_n^\l(s,X_s^n)\dif s+\sqrt2\int_0^t\nn u_n^\l(s,X_s^n)\dif W_s.\label{ito-krylov}
\end{align}
Hence, by the Burkholder-Davis-Gundy inequality, for $0\leq\tau_0\leq\tau_1\leq T$ with $\tau_1-\tau_0\leq\de$,
\begin{align}\label{basic-krylov}
    \Big\|\int_{\tau_0}^{\tau_1}f(s,X_s^n)\dif s\Big\|_{L^p(\Om)}
    \lesssim (1+\l\de)\|u_n^\l\|_\infty
    +\sqrt\de\|\nn u_n^\l\|_\infty.
\end{align}
Then \eqref{kol-star-Linfty}, \eqref{kol-star-grad} and \eqref{basic-krylov} give
\begin{align}
    \Big\|\int_{\tau_0}^{\tau_1}f(s,X_s^n)\dif s\Big\|_{L^p(\Om)}
    \lesssim
    \Big((1+\l\de)\ff{\l^{-1/2}}{w_{m,\eps,1}(J_\l)}
    +\ff{\sqrt\de}{\ell_{m+1}^{\eps}(J_\l)\ell_{m+2}(J_\l)}\Big)
    \|f\|_{\mL_T^\infty\bC^{-1,*}_{m,\eps}}.\no
\end{align}
Taking $\l=\l_0\vee\de^{-1}$ and using \eqref{Jl-asymp}, we get \eqref{krylov-star}.

Similarly, \eqref{kol-starstar-Linfty}, \eqref{kol-starstar-grad} and \eqref{basic-krylov} yield
\begin{align}
    \Big\|\int_{\tau_0}^{\tau_1}f(s,X_s^n)\dif s\Big\|_{L^p(\Om)}
    \lesssim
    \Big((1+\l\de)\ff{\l^{-1/2}}{w_{m,\eps,0}(J_\l)}
    +\ff{\sqrt\de}{\ell_{m+1}^{\eps}(J_\l)}\Big)
    \|f\|_{\mL_T^\infty\bC^{-1,**}_{m,\eps}}.\no
\end{align}
Again taking $\l=\l_0\vee\de^{-1}$ gives \eqref{krylov-starstar}. By a standard approximation, the above estimate still holds for
$f\in\mL_T^1C_b\cap\mL_T^\infty\bC^{-1,*}_{m,\eps}$ or
$f\in\mL_T^1 C_b\cap\mL_T^\infty\bC^{-1,**}_{m,\eps}$. The proof is complete.
\epf

\br
The dominant contribution in \eqref{krylov-star} and \eqref{krylov-starstar} comes from the martingale term in \eqref{ito-krylov}. More precisely,
\begin{align}
    (1+\l\de)\ff{\l^{-1/2}}{w_{m,\eps,1}(J_\l)}\bigg|_{\l=\de^{-1}}
    \asymp \bar\om_*(\de)\ll \om_*(\de),\ 
    (1+\l\de)\ff{\l^{-1/2}}{w_{m,\eps,0}(J_\l)}\bigg|_{\l=\de^{-1}}
    \asymp \bar\om_{**}(\de)\ll \om_{**}(\de).\no
\end{align}
\er

\bc\label{log-Af}
Let $p>2$. If $f\in\mL_T^\infty\bC^{-1,**}_{m,\eps}$ and there exist
$f_k\in\mL_T^1C_b\cap\mL_T^\infty\bC^{-1,**}_{m,\eps}$ such that
\begin{align}
    \lim_{k\to\infty}\|f_k-f\|_{\mL_T^\infty\bC^{-1,**}_{m,\eps}}=0,\no
\end{align}
then, for any fixed $n\in\mN$, the limit
\begin{align}
    (A_t^f)_n:=\lim_{k\to\infty}\int_0^t f_k(s,X_s^n)\dif s,\ t\in[0,T],\no
\end{align}
exists in $L^p(\Om;C([0,T]))$ and is independent of the choice of $\{f_k\}_{k\in\mN}$. Moreover, there exists a constant
$C:=C(p,T,d,m,\eps,\|b\|_{\mL_T^\infty\bC^{-1,*}_{m,\eps}})>0$ such that
\begin{align}\label{Af-starstar-est}
    \sup_{n\in\mN}\sup_{t\in[0,T]}\|(A_t^f)_n\|_{L^p(\Om)}
    +\sup_{n\in\mN}\sup_{\substack{s\neq t\in[0,T]}}
    \ff{\|(A_t^f)_n-(A_s^f)_n\|_{L^p(\Om)}}{\om_{**}(|t-s|)}
    \lesssim_{C}\|f\|_{\mL_T^\infty\bC^{-1,**}_{m,\eps}}.
\end{align}
If $f\in\mL_T^\infty\bC^{-1,*}_{m,\eps}$ and is approximable by functions in
$\mL_T^1C_b\cap\mL_T^\infty\bC^{-1,*}_{m,\eps}$ in the corresponding norm, then the term
$\om_{**}(|t-s|)\|f\|_{\mL_T^\infty\bC^{-1,**}_{m,\eps}}$ in \eqref{Af-starstar-est} can be replaced by
$\om_{*}(|t-s|)\|f\|_{\mL_T^\infty\bC^{-1,*}_{m,\eps}}$.
\ec

\bpf
Choose $q>p\vee2$. By the logarithmic Krylov estimate, for all $0\le s<t\le T$,
\begin{align}
\|A^{f_k-f_\ell}_t-A^{f_k-f_\ell}_s\|_{L^q(\Omega)}\lesssim \omega_{**}(t-s)\|f_k-f_\ell\|_{\mL_T^\infty\bC^{-1,**}_{m,\eps}}\lesssim |t-s|^{1/2}\|f_k-f_\ell\|_{\mL_T^\infty\bC^{-1,**}_{m,\eps}}. \nonumber
\end{align}
Hence the Kolmogorov criterion yields, for every $\gamma<\frac12-\frac1q$,
\begin{align}
 \|A^{f_k-f_\ell}\|_{L^q(\Omega;C^\gamma([0,T]))}
 \lesssim \|f_k-f_\ell\|_{\mL_T^\infty\bC^{-1,**}_{m,\eps}}. \nonumber
\end{align}
Since $\{f_k\}_{k\ge1}$ is Cauchy in $\mL_T^\infty\bC^{-1,**}_{m,\eps}$, $\{A^{f_k}\}_{k\ge1}$
is Cauchy in $L^q(\Omega;C([0,T]))$, hence also in $L^p(\Omega;C([0,T]))$.
\epf

We next prove the conditional version of Proposition \ref{log-krylov}.

\bl\label{cond-log-krylov}
Let $p\geq2$ and $0\leq\tau_0\leq\tau_1\leq T$ be two stopping times with $\tau_1-\tau_0\leq\de$. If
$f\in\mL_T^1C_b\cap\mL_T^\infty\bC^{-1,*}_{m,\eps}$, then
\begin{align}\label{cond-krylov-star}
    \sup_n\Big\|\mE\Big(\int_{\tau_0}^{\tau_1}f(s,X_s^n)\dif s\,\mid\,\cF_{\tau_0}\Big)\Big\|_{L^p(\Om)}
    \lesssim \bar\om_*(\de)\|f\|_{\mL_T^\infty\bC^{-1,*}_{m,\eps}}.
\end{align}
If $f\in\mL_T^1C_b\cap\mL_T^\infty\bC^{-1,**}_{m,\eps}$, then
\begin{align}
    \sup_n\Big\|\mE\Big(\int_{\tau_0}^{\tau_1}f(s,X_s^n)\dif s\,\mid\,\cF_{\tau_0}\Big)\Big\|_{L^p(\Om)}
    \lesssim \bar\om_{**}(\de)\|f\|_{\mL_T^\infty\bC^{-1,**}_{m,\eps}}.\no
\end{align}
Moreover, if $f$ belongs to either of the two spaces above and  can be approximated in the corresponding norm by functions in $\mL_T^1C_b$, then the same estimates hold with $\int_{\tau_0}^{\tau_1}f(s,X_s^n)\dif s$ replaced by $(A_{\tau_1}^f)_n-(A_{\tau_0}^f)_n$.
\el

\bpf
We only prove the estimate for $\bC^{-1,*}_{m,\eps}$, since the proof for $\bC^{-1,**}_{m,\eps}$ is similar. First assume that $f\in C_c^\infty([0,T]\times\mR^d)$. Let $u_n^\l$ solve the following backward equation.
\begin{align}
    \p_t u_n^\l+\Delta u_n^\l+b_n\cdot\nn u_n^\l-\l u_n^\l+f=0,
    \ u_{n}^\l(T)=0.\no
\end{align}
Applying the It\^o formula gives
\begin{align}
    \int_{\tau_0}^{\tau_1}f(s,X_s^n)\dif s
    &=u_n^\l(\tau_0,X_{\tau_0}^n)-u_n^\l(\tau_1,X_{\tau_1}^n)
    +\l\int_{\tau_0}^{\tau_1}u_n^\l(s,X_s^n)\dif s
    +\sqrt2\int_{\tau_0}^{\tau_1}\nn u_n^\l(s,X_s^n)\dif W_s.\no
\end{align}
Taking conditional expectation with respect to $\cF_{\tau_0}$ kills the martingale term. Hence
\begin{align}
    \Big\|\mE\Big(\int_{\tau_0}^{\tau_1}f(s,X_s^n)\dif s\,\big|\,\cF_{\tau_0}\Big)\Big\|_{L^p(\Om)}
    \lesssim (1+\l\de)\|u_n^\l\|_\infty.\no
\end{align}
By \eqref{kol-star-Linfty},
\begin{align}
    \Big\|\mE\Big(\int_{\tau_0}^{\tau_1}f(s,X_s^n)\dif s\,\big|\,\cF_{\tau_0}\Big)\Big\|_{L^p(\Om)}
    \lesssim (1+\l\de)\ff{\l^{-1/2}}{w_{m,\eps,1}(J_\l)}
    \|f\|_{\mL_T^\infty\bC^{-1,*}_{m,\eps}}.\no
\end{align}
Taking $\l=\l_0\vee\de^{-1}$ and using \eqref{Jl-asymp}, we obtain \eqref{cond-krylov-star}. Then by a standard approximation argument, the proof is complete.
\epf

The following proposition is crucial for Zvonkin transformation.
\bp\label{log-stoch-substitution}
Let $p>2$ and let $B$ be either $\mL_T^\infty\bC^{-1,*}_{m,\eps}$ or $\mL_T^\infty\bC^{-1,**}_{m,\eps}$. Let $f\in B$ be approximable by functions in $\mL_T^1C_b\cap B$. Let $g:[0,T]\times\mR^d\to\mR$ be bounded and satisfy
\begin{align}
    |g(t,x)-g(s,y)|\lesssim\big(|t-s|^{1/2}+|x-y|\big),\ s,t\in[0,T],\ x,y\in\mR^d.\no
\end{align}
For every $n\in\mN$ and $t\in[0,T]$, set
\begin{align}
    S_\pi^n(t):=\sum_{[t_i,t_{i+1}]\in\pi}g(t_i,X_{t_i}^n)\big((A_{t_{i+1}}^f)_n-(A_{t_i}^f)_n\big),\no
\end{align}
where $\pi$ is a partition of $[0,t]$. Then there exists an $L^p(\Om)$-limit, denoted by $\int_0^t g(s,X_s^n)\dif (A_s^f)_n$, uniformly in $n$ in the sense that
\begin{align}
    \lim_{|\pi|\to0}\sup_n\Big\|S_\pi^n(t)-\int_0^t g(s,X_s^n)\dif (A_s^f)_n\Big\|_{L^p(\Om)}=0.\no
\end{align}
Moreover,
\begin{align}\label{stoch-sub-formula}
    \int_0^t g(s,X_s^n)\dif (A_s^f)_n=(A_t^{gf})_n,
    \ t\in[0,T].
\end{align}
\ep

To prove this proposition, we need the following logarithmic stochastic sewing.

\bl\label{log-stoch-sewing}
Let $(\Om,\cF,(\cF)_{t\in[0,T]},\mP)$ be a filtered probability space, $p\geq2$ and let $B$ be either $\mL_T^\infty\bC^{-1,*}_{m,\eps}$ or $\mL_T^\infty\bC^{-1,**}_{m,\eps}$. Define
\begin{align}
    (\om_B,\bar\om_B):=\begin{cases}
    (\om_*,\bar\om_*),& B=\mL_T^\infty\bC^{-1,*}_{m,\eps},\\
    (\om_{**},\bar\om_{**}),& B=\mL_T^\infty\bC^{-1,**}_{m,\eps}.
    \end{cases}\no
\end{align}
Let $(\Xi_{s,t}^n)_{0\leq s\leq t\leq T,n\in\mN}$ be a family of real-valued random variables such that $\Xi_{s,s}^n=0$, $\Xi_{s,t}^n\in L^p(\Om)$, and $\Xi_{s,t}^n$ is $\cF_t$-measurable. Suppose that there exists a constant $C>0$, independent of $n$, such that for all $0\leq s\leq u\leq t\leq T$,
\begin{align}\label{log-sewing-ass1}
    \sup_n\|\de\Xi_{s,u,t}^n\|_{L^p(\Om)}
    \leq C |u-s|^{1/2}\om_B(t-u),
\end{align}
\begin{align}\label{log-sewing-ass2}
    \sup_n\big\|\mE\big(\de\Xi_{s,u,t}^n\,\big|\,\cF_s\big)\big\|_{L^p(\Om)}
    \leq C |u-s|^{1/2}\bar\om_B(t-u),
\end{align}
where
\begin{align}
    \de\Xi_{s,u,t}^n:=\Xi_{s,t}^n-\Xi_{s,u}^n-\Xi_{u,t}^n.\no
\end{align}
Then, for every $t\in[0,T]$, there exist $\cF_t$-measurable random variables $\cI_t^n(\Xi)\in L^p(\Om)$ such that, for any deterministic partition $\pi$ of $[0,t]$,
\begin{align}\label{log-sewing-limit}
\lim_{|\pi|\to0}\sup_n\Big\|\cI_t^n(\Xi)-\sum_{[t_i,t_{i+1}]\in\pi}\Xi_{t_i,t_{i+1}}^n\Big\|_{L^p(\Om)}=0.
\end{align}
Moreover, the map $\Xi\mapsto\cI^n(\Xi)$ is linear and satisfies
\begin{align}\label{log-sewing-cont}
    \sup_n\|\cI_t^n(\Xi)-\Xi_{0,t}^n\|_{L^p(\Om)}\lesssim_{p,T} C.
\end{align}
\el

\bpf
This is the usual stochastic sewing argument of \cite{Le20}, with the geometric summability replaced by the logarithmic summability below; see also \cite{MP24} for related extensions. Set
\begin{align}
    \eta_B(r):=r^{1/2}\om_B(r),\ \bar\eta_B(r):=r^{1/2}\bar\om_B(r),\ r>0,\no
\end{align}
and
\begin{align}
    \ell_B(r):=\begin{cases}
    \ell_{m+2}^{\eps}(r^{-1})\ell_{m+3}(r^{-1}),& B=\mL_T^\infty\bC^{-1,*}_{m,\eps},\\
    \ell_{m+2}^{\eps}(r^{-1}),& B=\mL_T^\infty\bC^{-1,**}_{m,\eps}.
    \end{cases}\no
\end{align}
A direct calculation gives, uniformly for $0<r\leq T$,
\begin{align}
    \sum_{k\geq0}2^{k/2}\eta_B(2^{-k}r)
    +\sum_{k\geq0}2^k\bar\eta_B(2^{-k}r)
    \lesssim_T \ff{r}{\ell_B(r)}.\no
\end{align}
Let $\pi_k(s,t)$ be the dyadic partition of $[s,t]$ with mesh $h_k:=2^{-k}(t-s)$, and set
\begin{align}
    S_k^n(s,t):=\sum_{[r,v]\in\pi_k(s,t)}\Xi_{r,v}^n.\no
\end{align}
Decomposing the increment of $S_{k+1}^n(s,t)-S_k^n(s,t)$ into its martingale part and its conditional expectation part: if $r_i=s+ih_k$ and $q_i=(r_i+r_{i+1})/2$, then
\begin{align}
    S_{k+1}^n(s,t)-S_k^n(s,t)
    &=-\sum_{i=0}^{2^k-1}\de\Xi_{r_i,q_i,r_{i+1}}^n\no\\
    &=-\sum_{i=0}^{2^k-1}\Big[\big(\de\Xi_{r_i,q_i,r_{i+1}}^n-\mE(\de\Xi_{r_i,q_i,r_{i+1}}^n\mid\cF_{r_i})\big)+\mE(\de\Xi_{r_i,q_i,r_{i+1}}^n\mid\cF_{r_i})\Big].\no
\end{align}
The martingale part $\de\Xi_{r_i,q_i,r_{i+1}}^n-\mE(\de\Xi_{r_i,q_i,r_{i+1}}^n\mid\cF_{r_i})$ are adapted to the filtration $(\cF_{r_{i+1}})_i$ with
\begin{align}
\|\de\Xi_{r_i,q_i,r_{i+1}}^n-\mE(\de\Xi_{r_i,q_i,r_{i+1}}^n\mid\cF_{r_i})\|_{L^p(\Om)}
\leq2\|\de\Xi_{r_i,q_i,r_{i+1}}^n\|_{L^p(\Om)}
\overset{\eqref{log-sewing-ass1}}{\lesssim} C\eta_B(h_k).\no
\end{align}
Hence, by the Burkholder-Davis-Gundy inequality and then the Minkowski inequality,
\begin{align}
&\Big\|\sum_{i=0}^{2^k-1}\big(\de\Xi_{r_i,q_i,r_{i+1}}^n-\mE(\de\Xi_{r_i,q_i,r_{i+1}}^n\mid\cF_{r_i})\big)\Big\|_{L^p(\Om)}\no\\
&\ \lesssim
\Big\|\Big(\sum_{i=0}^{2^k-1}\big|\de\Xi_{r_i,q_i,r_{i+1}}^n-\mE(\de\Xi_{r_i,q_i,r_{i+1}}^n\mid\cF_{r_i})\big|^2\Big)^{1/2}\Big\|_{L^p(\Om)}\no\\
&\ \lesssim
\Big(\sum_{i=0}^{2^k-1}\|\de\Xi_{r_i,q_i,r_{i+1}}^n-\mE(\de\Xi_{r_i,q_i,r_{i+1}}^n\mid\cF_{r_i})\|_{L^p(\Om)}^2\Big)^{1/2}
\lesssim 2^{k/2}C\eta_B(h_k).\no
\end{align}
As for the conditional expectation part, similarly by \eqref{log-sewing-ass2}, we obtain
\begin{align}
\Big\|\sum_{i=0}^{2^k-1}\mE(\de\Xi_{r_i,q_i,r_{i+1}}^n\mid\cF_{r_i})\Big\|_{L^p(\Om)}
&\lesssim\sum_{i=0}^{2^k-1}\|\mE(\de\Xi_{r_i,q_i,r_{i+1}}^n\mid\cF_{r_i})\|_{L^p(\Om)}\no\\
&\lesssim\sum_{i=0}^{2^k-1}\bar\eta_B(h_k)\lesssim 2^kC\bar\eta_B(h_k).\no
\end{align}
These yield
\begin{align}\label{dyadic-sewing-est}
    \sup_n\|S_{k+1}^n(s,t)-S_k^n(s,t)\|_{L^p(\Om)}
    \lesssim C\Big(2^{k/2}\eta_B(h_k)+2^k\bar\eta_B(h_k)\Big).
\end{align}

For $N<M$, by \eqref{dyadic-sewing-est},
\begin{align}\label{dyadic-Cauchy-tail}
    \sup_n\|S_M^n(s,t)-S_N^n(s,t)\|_{L^p(\Om)}
    \leq \sum_{k=N}^{M-1}\sup_n\|S_{k+1}^n(s,t)-S_k^n(s,t)\|_{L^p(\Om)}
    \lesssim C R_N(t-s),
\end{align}
where
\begin{align}
    R_N(r):=\sum_{k=N}^{\infty}
    \Big(2^{k/2}\eta_B(2^{-k}r)+2^k\bar\eta_B(2^{-k}r)\Big),
    \ 0<r\leq T.\no
\end{align}
Noting $R_N(r)$ is the tail of a convergent positive series, we have for $0<r\leq T$,
\begin{align}
    R_N(r)\rightarrow0,\ N\to\infty.\no
\end{align}
Hence \eqref{dyadic-Cauchy-tail} implies
\begin{align}
    \lim_{N,M\to\infty}
    \sup_n\|S_M^n(s,t)-S_N^n(s,t)\|_{L^p(\Om)}=0.\no
\end{align}
Thus $(S_k^n(s,t))_{k\geq0}$ is a Cauchy sequence in $L^p(\Om)$, uniformly in $n$. Since $L^p(\Om)$ is complete, for each $n$ there exists an $L^p$-limit, denoted by
\begin{align}
    \cI_{s,t}^n(\Xi):=\lim_{k\to\infty}S_k^n(s,t)
    \ \hbox{in }L^p(\Om).\no
\end{align}
Moreover, letting $M\to\infty$ in \eqref{dyadic-Cauchy-tail} gives
\begin{align}
    \sup_n\|\cI_{s,t}^n(\Xi)-S_N^n(s,t)\|_{L^p(\Om)}
    \lesssim_p C R_N(t-s).\no
\end{align}

Denote $\cI_t^n(\Xi):=\cI_{0,t}^n(\Xi)$. By the deterministic partition decomposition in \cite[Lemma 2.1]{MP24}, the same martingale-difference and conditional-expectation estimates as in \eqref{dyadic-sewing-est}, together with the logarithmic summability above, yield, for every deterministic finite partition $\pi$ of $[0,t]$,
\begin{align}\label{general-partition-sewing-est}
\sup_n\Big\|\cI_t^n(\Xi)-\sum_{[r,v]\in\pi}\Xi_{r,v}^n\Big\|_{L^p(\Om)}\lesssim_{p,T}C\sum_{[r,v]\in\pi}\ff{v-r}{\ell_B(v-r)}\leq\ff{C_{p,T}Ct}{\ell_B(|\pi|)}.
\end{align}
Since $\ell_B(r)\to\infty$ as $r\downarrow0$, this proves \eqref{log-sewing-limit}. Taking the one-interval partition in \eqref{general-partition-sewing-est} gives \eqref{log-sewing-cont}. Finally, $\cI_t^n(\Xi)$ is $\cF_t$-measurable and depends linearly on $\Xi$, since it is the $L^p(\Om)$-limit of the dyadic sums $S_k^n(0,t)$.
\epf

Now we return to the proof of Proposition \ref{log-stoch-substitution}.

\bpf[Proof of Proposition \ref{log-stoch-substitution}]
We first note that multiplication by $g$ is bounded on the spaces $\bC^{-1,*}_{m,\eps}$ and $\bC^{-1,**}_{m,\eps}$, uniformly in time. More precisely, by the usual paraproduct multiplier estimate,
\begin{align}\label{multiplier-g}
    \|gf\|_B\lesssim_g\|f\|_B,
\end{align}
and if $f_k\to f$ in $B$, then $gf_k\to gf$ in $B$. Indeed, by the Bony decomposition,
$$
gf=g\prec f+g\circ f+f\prec g .
$$
For $g\prec f$, the estimate follows directly from $\|S_{j-1}g\|_\infty\leq\|g\|_\infty$ and the finite overlap of dyadic blocks. For the remaining two terms, the Bernstein inequality and $\|\cR_jg\|_\infty\lesssim2^{-j}\|\nn g\|_\infty$ give
$$
w_{m,\eps,\theta}(j)2^{-j}
\|\cR_j(g\circ f+f\prec g)\|_\infty
\lesssim
\|\nn g\|_\infty
\|f\|_{\bC^{-1}_{m,\eps,\theta}},
\ \theta\in\{0,1\}.
$$
Taking the supremum over $j$ and then over $t\in[0,T]$ proves the claim.

For $0\leq s\leq t\leq T$, define
\begin{align}
    \Xi_{s,t}^n:=g(s,X_s^n)\big((A_t^f)_n-(A_s^f)_n\big).\no
\end{align}
Then $\Xi_{s,s}^n=0$, $\Xi_{s,t}^n$ is $\cF_t$-measurable, and for $0\leq s\leq u\leq t\leq T$,
\begin{align}\label{delta-Xi-sub}
    \de\Xi_{s,u,t}^n
    =\big(g(s,X_s^n)-g(u,X_u^n)\big)\big((A_t^f)_n-(A_u^f)_n\big).
\end{align}
By Proposition \ref{log-krylov},
\begin{align}
    \sup_n\|X_t^n-X_s^n\|_{L^{2p}(\Om)}
    \lesssim |t-s|^{1/2},\ 0\leq s\leq t\leq T.\no
\end{align}
Therefore the assumption on $g$ gives
\begin{align}\label{g-holder-moment}
    \sup_n\|g(t,X_t^n)-g(s,X_s^n)\|_{L^{2p}(\Om)}
    \lesssim |t-s|^{1/2}.
\end{align}
On the other hand, Corollary \ref{log-Af} gives
\begin{align}\label{Af-increment-B}
    \sup_n\|(A_t^f)_n-(A_s^f)_n\|_{L^{2p}(\Om)}
    \lesssim \om_B(t-s)\|f\|_B.
\end{align}
Combining \eqref{delta-Xi-sub}, \eqref{g-holder-moment} and \eqref{Af-increment-B}, we get
\begin{align}
    \sup_n\|\de\Xi_{s,u,t}^n\|_{L^p(\Om)}
    \lesssim |u-s|^{1/2}\om_B(t-u)\|f\|_B.\no
\end{align}

We next estimate the conditional expectation. Since $g(s,X_s^n)-g(u,X_u^n)$ is $\cF_u$-measurable,
\begin{align}
    \mE(\de\Xi_{s,u,t}^n\mid\cF_s)
    =\mE\Big((g(s,X_s^n)-g(u,X_u^n))
    \mE\big((A_t^f)_n-(A_u^f)_n\mid\cF_u\big)\,\Big|\,\cF_s\Big).\no
\end{align}
Using the H\"older inequality, \eqref{g-holder-moment}, and Lemma \ref{cond-log-krylov}, we obtain
\begin{align}
    \sup_n\big\|\mE(\de\Xi_{s,u,t}^n\mid\cF_s)\big\|_{L^p(\Om)}
    \lesssim |u-s|^{1/2}\bar\om_B(t-u)\|f\|_B.\no
\end{align}
Thus the assumptions of Lemma \ref{log-stoch-sewing} are satisfied with parameter $C\lesssim\|f\|_B$. Consequently, for every $t\in[0,T]$, the Riemann sums $S_\pi^n(t)$ converge in $L^p(\Om)$ as $|\pi|\to0$, uniformly in $n$. This proves the existence of the left hand side of \eqref{stoch-sub-formula}. Moreover, by \eqref{log-sewing-cont}, \eqref{Af-increment-B}, and the linearity of the sewing map, for every approximable $h\in B$,
\begin{align}\label{sewn-cont-in-f}
    \sup_n\Big\|\int_0^t g(s,X_s^n)\dif(A_s^h)_n\Big\|_{L^p(\Om)}
    \lesssim_g \|h\|_B.
\end{align}

It remains to identify the limit. First assume that $f\in\mL_T^1C_b\cap B$. Then $(A_t^f)_n=\int_0^t f(r,X_r^n)\dif r$, and for every partition $\pi$ of $[0,t]$,
\begin{align}
    S_\pi^n(t)
    =\sum_{[t_i,t_{i+1}]\in\pi}g(t_i,X_{t_i}^n)\int_{t_i}^{t_{i+1}}f(r,X_r^n)\dif r.\no
\end{align}
Since $X^n$ has continuous paths and $g$ is continuous, the left-point step process
\begin{align}
    r\mapsto \sum_{[t_i,t_{i+1}]\in\pi}g(t_i,X_{t_i}^n)\bbone_{[t_i,t_{i+1})}(r)\no
\end{align}
converges uniformly on $[0,t]$, almost surely, to $r\mapsto g(r,X_r^n)$. Since $r\mapsto\|f(r)\|_\infty$ belongs to $L^1([0,T])$, it follows that
\begin{align}
    S_\pi^n(t)\rightarrow \int_0^t g(r,X_r^n)f(r,X_r^n)\dif r
    =(A_t^{gf})_n\no
\end{align}
in probability. The same sums converge in $L^p(\Om)$ to the sewn limit; hence the sewn limit must be $(A_t^{gf})_n$, and \eqref{stoch-sub-formula} holds for $f\in\mL_T^1C_b\cap B$.

For general approximable $f\in B$, take $f_k\in\mL_T^1C_b\cap B$ with $\|f_k-f\|_B\to0$. The identity already proved gives
\begin{align}
    \int_0^t g(s,X_s^n)\dif(A_s^{f_k})_n=(A_t^{gf_k})_n.\no
\end{align}
By \eqref{sewn-cont-in-f},
\begin{align}
    \lim_{k\to\infty}\sup_n\Big\|\int_0^t g(s,X_s^n)\dif(A_s^{f_k-f})_n\Big\|_{L^p(\Om)}=0.\no
\end{align}
By Corollary \ref{log-Af} and \eqref{multiplier-g},
\begin{align}
    \lim_{k\to\infty}\sup_n\|(A_t^{g(f_k-f)})_n\|_{L^p(\Om)}=0.\no
\end{align}
Letting $k\to\infty$ in the identity for $f_k$ yields \eqref{stoch-sub-formula} for $f$. The proof is complete.
\epf

\section{Solvability of SDEs with distributional drifts}\label{sec5}

We now prove the main theorem. The proof is divided into three parts. First, we construct weak solutions by compactness from the smooth approximating equations, identify the limiting drift through the logarithmic Krylov estimates, and show that the constructed solutions belong to the logarithmic Krylov class. Second, we prove uniqueness in law within this class by means of the Zvonkin transformation. Finally, we establish the two-sided Gaussian bounds by applying the Aronson estimates to the smooth divergence-form approximations and passing to the limit.

\subsection{Existence of weak solutions}

First, we prove the existence of the solutions to SDE \eqref{SDE} by a standard weak convergence argument.

\bt\label{main-existence}
Assume that $\eps>0$, $m\in\mN$, $b\in\mL_T^\infty\bC^{-1,*}_{m,\eps}$ and $\div b=0$. For any $\mu\in\sP(\mR^d)$, there exists a weak solution $(\mathfrak F,X,W)$ to the SDE \eqref{SDE} in the sense of Definition \ref{weaksol}.
\et

\bpf
Take a probability space $(\Om^0,\sF^0,(\sF_t^0)_{t\geq0},\bP^0)$ carrying an $\mR^d$-valued random variable $X_0$ with law $\mu$ and a standard $d$-dimensional Brownian motion $W^0$ independent of $X_0$. For every $n\in\mN$, let $X^n$ be the unique strong solution to
\begin{align}
    X_t^n=X_0+\int_0^t b_n(s,X_s^n)\dif s+\sqrt2 W_t^0,
    \  t\in[0,T].\no
\end{align}
Set $A_t^n:=\int_0^t b_n(s,X_s^n)\dif s$. By Proposition \ref{log-krylov} and the uniform bound
\begin{align}
    \sup_n\|b_n\|_{\mL_T^\infty\bC^{-1,*}_{m,\eps}}
    \lesssim \|b\|_{\mL_T^\infty\bC^{-1,*}_{m,\eps}},\no
\end{align}
for every $p>2$,
\begin{align}
    \sup_n\|A_t^n-A_s^n\|_{L^p(\Om^0)}
    \lesssim \om_*(|t-s|)\|b\|_{\mL_T^\infty\bC^{-1,*}_{m,\eps}},
    \ s,t\in[0,T].\no
\end{align}
Hence, by the Kolmogorov continuity criterion, for every $\eta\in(0,1/2-1/p)$,
\begin{align}
    \sup_n\mE\|A^n\|_{C^\eta([0,T])}^p<\infty,
    \ \mE\|W^0\|_{C^\eta([0,T])}^p<\infty.\no
\end{align}
Using also the tightness of $\mu$, we obtain the tightness of the laws of $\{X^n\}_{n\in\mN}$.

Define
\begin{align}
    Q^n:=\bP^0\circ(X^n,W^0)^{-1}\in\sP(C_T\times C_T).\no
\end{align}
Since the marginal law of $W^0$ is fixed and tight, $\{Q^n\}_{n\in\mN}$ is tight on $C_T\times C_T$. By the Prohorov theorem and Skorokhod representation theorem, there exist a probability space $(\Om,\sF,\bP)$ and $C_T\times C_T$-valued random variables $(X^n,W^n)$ and $(X,W)$ such that, along a subsequence,
\begin{align}
    \bP\circ(X^n,W^n)^{-1}=Q^n,
    \ \bP\circ(X,W)^{-1}=Q,\no
\end{align}
for some $Q\in\sP(C_T\times C_T)$, and
\begin{align}
    (X^n,W^n)\to(X,W)\ \hbox{in }C_T\times C_T,
    \ \bP\hbox{-a.s.}\no
\end{align}
In particular,
\begin{align}
    X_t^n=X_0^n+\int_0^t b_n(s,X_s^n)\dif s+
    \sqrt2 W_t^n,
    \ t\in[0,T],\ \bP\hbox{-a.s.}\no
\end{align}
Moreover, $\bP\circ X_0^{-1}=\mu$. Let $\sF_t$ be the usual augmentation of the natural filtration generated by $(X,W)$. Since each $W^n$ is a Brownian motion with respect to the natural filtration of $(X^n,W^n)$ and the Brownian martingale problem is stable under weak convergence, $W$ is an $\sF_t$-Brownian motion.

Next we identify the limiting drift. First note that
\begin{align}\label{bn-to-b-starstar}
    \lim_{n\to\infty}\|b_n-b\|_{\mL_T^\infty\bC^{-1,**}_{m,\eps}}=0.
\end{align}
Indeed, for the high frequency tail,
\begin{align}
    \sup_{j\geq N}w_{m,\eps,0}(j)2^{-j}\|\cR_j b\|_{\mL_T^\infty L^\infty}
    \leq \ff1{\ell_{m+2}(N)}\|b\|_{\mL_T^\infty\bC^{-1,*}_{m,\eps}}\to0,
    \ N\to\infty,\no
\end{align}
and the finite number of dyadic blocks converges directly by mollification.

Let $f\in\mL_T^1C_b\cap\mL_T^\infty\bC^{-1,**}_{m,\eps}$. By dominated convergence,
\begin{align}
    \int_s^t f(r,X_r^n)\dif r
    \to
    \int_s^t f(r,X_r)\dif r
    \ \hbox{in }L^p(\Om),\ 0\leq s\leq t\leq T.\no
\end{align}
Passing to the limit in Proposition \ref{log-krylov}, we get
\begin{align}
    \Big\|\int_s^t f(r,X_r)\dif r\Big\|_{L^p(\Om)}
    \lesssim
    \om_{**}(t-s)\|f\|_{\mL_T^\infty\bC^{-1,**}_{m,\eps}},
    \ 0\leq s\leq t\leq T.\no
\end{align}
Consequently, if $f_k\to f$ in $\mL_T^\infty\bC^{-1,**}_{m,\eps}$ with $f_k\in\mL_T^1C_b\cap\mL_T^\infty\bC^{-1,**}_{m,\eps}$, then $\int_0^t f_k(s,X_s)\dif s$ is Cauchy in $L^2(\Om)$, and the limit is independent of the approximating sequence. In particular, by \eqref{bn-to-b-starstar}, we may define
\begin{align}
    A_t^b:=\lim_{k\to\infty}\int_0^t b_k(s,X_s)\dif s
    \ \hbox{in }L^2(\Om).\no
\end{align}
The same estimate with $p>2$ and the Kolmogorov criterion gives a continuous version of $t\mapsto A_t^b$.

It remains to prove that the drift term in the approximating equation converges to $A_t^b$. Fix $t\in[0,T]$ and $k\in\mN$. We decompose
\begin{align}
    &\Big\|\int_0^t b_n(s,X_s^n)\dif s-A_t^b\Big\|_{L^2(\Om)}\no\\
    &\ \leq
    \Big\|\int_0^t (b_n-b_k)(s,X_s^n)\dif s\Big\|_{L^2(\Om)}
    +\Big\|\int_0^t b_k(s,X_s^n)\dif s-
    \int_0^t b_k(s,X_s)\dif s\Big\|_{L^2(\Om)}\no\\
    &\ 
    +\Big\|\int_0^t b_k(s,X_s)\dif s-A_t^b\Big\|_{L^2(\Om)}.\no
\end{align}
For the first term, Proposition \ref{log-krylov} gives
\begin{align}
    \limsup_{n\to\infty}
    \Big\|\int_0^t (b_n-b_k)(s,X_s^n)\dif s\Big\|_{L^2(\Om)}
    \lesssim
    \|b-b_k\|_{\mL_T^\infty\bC^{-1,**}_{m,\eps}}.\no
\end{align}
For the second term, since $b_k\in\mL_T^\infty C_b^\infty$ and $X^n\to X$ almost surely in $C_T$, the dominated convergence theorem yields
\begin{align}
    \lim_{n\to\infty}
    \Big\|\int_0^t b_k(s,X_s^n)\dif s-
    \int_0^t b_k(s,X_s)\dif s\Big\|_{L^2(\Om)}=0.\no
\end{align}
The third term tends to zero as $k\to\infty$ by the definition of $A_t^b$. Therefore, first taking $\limsup_{n\to\infty}$ and then letting $k\to\infty$, we obtain
\begin{align}
    \lim_{n\to\infty}
    \Big\|\int_0^t b_n(s,X_s^n)\dif s-A_t^b\Big\|_{L^2(\Om)}=0.\no
\end{align}
Passing to the limit in the approximating equation, we get
\begin{align}
    X_t=X_0+A_t^b+\sqrt2 W_t,
    \ t\in[0,T],\ \bP\hbox{-a.s.}\no
\end{align}
Thus $(\mathfrak F,X,W)$, with $\mathfrak F:=(\Om,\sF,(\sF_t)_{t\geq0},\bP)$, is a weak solution to SDE \eqref{SDE}. The proof is complete.
\epf

We next isolate the Krylov class used in the uniqueness part.

\bd\label{krylov-class-def}
A weak solution $(\mathfrak F,X,W)$ to SDE \eqref{SDE} is said to satisfy the logarithmic Krylov estimate if for every $p>2$ and
$f\in\mL_T^1C_b\cap\mL_T^\infty\bC^{-1,**}_{m,\eps}$, there exists a constant
$C=C(p,T,d,m,\eps,\|b\|_{\mL_T^\infty\bC^{-1,*}_{m,\eps}})>0$ such that, for all
$0\leq s\leq t\leq T$,
\begin{align}\label{weak-krylov-est}
    \Big\|\int_s^t f(r,X_r)\dif r\Big\|_{L^p(\Om)}
    \lesssim_C
    \om_{**}(t-s)\|f\|_{\mL_T^\infty\bC^{-1,**}_{m,\eps}},
\end{align}
and
\begin{align}\label{weak-cond-krylov-est}
    \Big\|\mE\Big(\int_s^t f(r,X_r)\dif r\,\big|\,\cF_s\Big)\Big\|_{L^p(\Om)}
    \lesssim_C
    \bar\om_{**}(t-s)\|f\|_{\mL_T^\infty\bC^{-1,**}_{m,\eps}}.
\end{align}
\ed

\bp\label{constructed-solution-krylov}
The weak solution constructed in Theorem \ref{main-existence} satisfies Definition \ref{krylov-class-def}.
\ep

\bpf
We keep the notation in the proof of Theorem \ref{main-existence}. First let
$f\in\mL_T^1C_b\cap\mL_T^\infty\bC^{-1,**}_{m,\eps}$. For the approximating solutions $X^n$, Proposition \ref{log-krylov} and Lemma \ref{cond-log-krylov} give, uniformly in $n$,
\begin{align}\label{approx-krylov-for-limit-proof}
    \Big\|\int_s^t f(r,X_r^n)\dif r\Big\|_{L^p(\Om)}
    \lesssim
    \om_{**}(t-s)\|f\|_{\mL_T^\infty\bC^{-1,**}_{m,\eps}},
\end{align}
and
\begin{align}\label{approx-cond-krylov-for-limit-proof}
    \Big\|\mE\Big(\int_s^t f(r,X_r^n)\dif r\,\big|\,\cF_s^n\Big)\Big\|_{L^p(\Om)}
    \lesssim
    \bar\om_{**}(t-s)\|f\|_{\mL_T^\infty\bC^{-1,**}_{m,\eps}},
\end{align}
where $\cF_s^n$ denotes the natural filtration generated by $(X^n,W^n)$.

Since $(X^n,W^n)\to(X,W)$ almost surely in $C([0,T];\mR^d)\times C([0,T];\mR^d)$, and since
$f\in\mL_T^1C_b$, the dominated convergence theorem yields
\begin{align}
    \int_s^t f(r,X_r^n)\dif r
    \to
    \int_s^t f(r,X_r)\dif r
    \ \hbox{in }L^p(\Om).\no
\end{align}
Passing to the limit in \eqref{approx-krylov-for-limit-proof}, we obtain \eqref{weak-krylov-est}.

It remains to pass the conditional estimate to the limit. Let $q=p/(p-1)$ and let $\xi$ be a bounded continuous $\cF_s$-measurable cylinder random variable. Then, by \eqref{approx-cond-krylov-for-limit-proof},
\begin{align}
    \bigg|\mE\bigg[\xi(X^n,W^n)\int_s^t f(r,X_r^n)\dif r\bigg]\bigg|
    &\leq
    \|\xi(X^n,W^n)\|_{L^q(\Om)}
    \Big\|\mE\Big(\int_s^t f(r,X_r^n)\dif r\,\big|\,\cF_s^n\Big)\Big\|_{L^p(\Om)}
    \no\\
    &\lesssim
    \|\xi(X^n,W^n)\|_{L^q(\Om)}
    \bar\om_{**}(t-s)\|f\|_{\mL_T^\infty\bC^{-1,**}_{m,\eps}}.\no
\end{align}
Letting $n\to\infty$ gives
\begin{align}
    \bigg|\mE\bigg[\xi(X,W)\int_s^t f(r,X_r)\dif r\bigg]\bigg|
    \lesssim
    \|\xi(X,W)\|_{L^q(\Om)}
    \bar\om_{**}(t-s)\|f\|_{\mL_T^\infty\bC^{-1,**}_{m,\eps}}.\no
\end{align}
By the density of bounded continuous cylinder random variables in $L^q(\cF_s)$ and duality, this implies
\begin{align}
    \Big\|\mE\Big(\int_s^t f(r,X_r)\dif r\,\big|\,\cF_s\Big)\Big\|_{L^p(\Om)}
    \lesssim
    \bar\om_{**}(t-s)\|f\|_{\mL_T^\infty\bC^{-1,**}_{m,\eps}}.\no
\end{align}
This proves \eqref{weak-cond-krylov-est}. The proof is complete.
\epf

\subsection{Uniqueness of the solution}

We now remove the distributional drift by a Zvonkin transform. To this end, we solve the backward Kolmogorov equation with source term $b$. By Corollary \ref{endpoint-kolmogorov}, there exists a unique solution $u^\l$ satisfying the backward equation
\begin{align}\label{zvonkin-pde}
    \p_tu^\l+\Delta u^\l-\l u^\l+b\cdot\nn u^\l+b=0,
    \  u^\l(T)=0,
\end{align}
where the term $b\cdot\nn u^\l$ is defined in the sense of Lemma \ref{driftterm}. Moreover,
\begin{align}
    \|\nn u^\l\|_\infty
    \lesssim
    \ff1{\ell_{m+1}^{\eps}(J_\l)\ell_{m+2}(J_\l)}
    \|b\|_{\mL_T^\infty\bC^{-1,*}_{m,\eps}}.\no
\end{align}
Choosing $\l$ large enough, we may assume
\begin{align}
    \|\nn u^\l\|_\infty\leq\ff12.\no
\end{align}
Let
\begin{align}
    \Phi(t,x):=x+u^\l(t,x).\no
\end{align}
Then for each $t\in[0,T]$, $x\mapsto\Phi(t,x)$ is a $C^1$-diffeomorphism and
\begin{align}
    \|\nn\Phi\|_\infty+\|\nn\Phi^{-1}\|_\infty\leq4.\no
\end{align}
We now show the crucial Zvonkin transformation.

\bl\label{zvonkin-transform}
Let $(\mathfrak F,X,W)$ be a weak solution of SDE \eqref{SDE} satisfying Definition \ref{krylov-class-def}. Then $Y_t:=\Phi(t,X_t)$ solves the following SDE:
\begin{align}
    \dif Y_t=\wt b(t,Y_t)\dif t+\wt\s(t,Y_t)\dif W_t,\no
\end{align}
where
\begin{align}
    \wt b(t,y):=\l u^\l(t,\Phi^{-1}(t,y)),\ 
    \wt \s(t,y):=\sqrt{2}\nn\Phi(t,\Phi^{-1}(t,y)).\no
\end{align}
Moreover, $\wt b$ and $\wt\s$ are bounded measurable, and for all $t\in[0,T]$, $y\neq y'\in\mR^d$,
\begin{align}
    |\wt b(t,y)-\wt b(t,y')|
    \leq C|y-y'|,\ 
    |\wt \s(t,y)-\wt \s(t,y')|
    \leq C\left(1\wedge
    \ff1{\ell_{m+2}^{\eps}(|y-y'|^{-1})}\right),\no
\end{align}
where $C:=C(T,d,m,\eps,\|b\|_{\mL_T^\infty\bC^{-1,*}_{m,\eps}})$, and for all $\xi\in\mR^d$,
\begin{align}
\frac{|\xi|}{4}\leq|\wt\s(t,y)\xi|\leq4|\xi|.\no
\end{align}
\el

\bpf
Since $X$ satisfies the logarithmic Krylov estimate in Definition \ref{krylov-class-def}, the arguments of Corollary \ref{log-Af} and Proposition \ref{log-stoch-substitution} apply to $X$ itself. In particular, for every approximable $f\in\mL_T^\infty\bC^{-1,**}_{m,\eps}$, we can define
\begin{align}
    A_t^f:=\lim_{k\to\infty}\int_0^t f_k(s,X_s)\dif s,
    \  t\in[0,T],\no
\end{align}
where $f_k\in\mL_T^1C_b\cap\mL_T^\infty\bC^{-1,**}_{m,\eps}$ and $f_k\to f$ in $\mL_T^\infty\bC^{-1,**}_{m,\eps}$. Moreover, the stochastic substitution formula holds:
\begin{align}\label{substitution-for-X}
    \int_0^t g(s,X_s)\dif A_s^f=A_t^{gf},
    \  t\in[0,T],
\end{align}
whenever $g$ is bounded and satisfies the regularity assumption in Proposition \ref{log-stoch-substitution}.

Let $u_n^\l(t,x):=u^\l(t,\cdot)*\phi_n(x)$ and $\Phi_n(t,x):=x+u_n^\l(t,x)$. By \eqref{zvonkin-pde}, after convolution in the space variable,
\begin{align}
    \p_t u_n^\l+\Delta u_n^\l-\l u_n^\l
    =-\phi_n*(b\cdot\nn u^\l+b).\no
\end{align}
Since $b\in\mL_T^\infty\bC^{-1,*}_{m,\eps}$ is approximable in the weaker space $\mL_T^\infty\bC^{-1,**}_{m,\eps}$, the additive functional $A^b$ is well defined along $X$ by the extended Krylov estimate. Moreover, for every partition $\pi$ of $[0,T]$,
$$
\mE\sum_{[s,t]\in\pi}|A_t^b-A_s^b|^2
\lesssim
\sum_{[s,t]\in\pi}\om_{**}^2(t-s)
\|b\|_{\mL_T^\infty\bC^{-1,**}_{m,\eps}}^2
\lesssim
\ff{T}{\ell_{m+2}^{2\eps}(|\pi|^{-1})}
\|b\|_{\mL_T^\infty\bC^{-1,**}_{m,\eps}}^2
\to0.
$$
Thus $A^b$ has zero quadratic variation, and $ X_t=X_0+A_t^b+\sqrt2 W_t$ is a Dirichlet process with quadratic variation $2tI_d$. For each fixed $n$, $\Phi_n$ is smooth in $x$ and absolutely continuous in $t$, with bounded derivatives. Hence, arguing as in \cite[Proof of Lemma 4.5]{HZ25} and applying the generalized It\^o formula (cf. \cite{Fol81}) to $\Phi_n(t,X_t)$, we obtain
\begin{align}
    \Phi_n(t,X_t)
    &=\Phi_n(0,X_0)+\l\int_0^t u_n^\l(s,X_s)\dif s
      +A_t^{r_n}
      +\sqrt2\int_0^t\nn\Phi_n(s,X_s)\dif W_s,\no
\end{align}
where $A_t^{r_n}$ is understood in the sense of \eqref{substitution-for-X} and 
\begin{align}
    r_n:=b\cdot\nn\Phi_n-\phi_n*(b\cdot\nn u^\l+b).\no
\end{align}
Noting $b\cdot\nn\Phi_n=b+b\cdot\nn u_n^\l$, we have
\begin{align}
    r_n
    =b-b*\phi_n+b\cdot\nn(u_n^\l-u^\l)
      +b\cdot\nn u^\l-\phi_n*(b\cdot\nn u^\l).\no
\end{align}
By Lemma \ref{driftterm}, $b\cdot\nn u^\l\in\mL_T^\infty\bC^{-1,*}_{m,\eps}$. We claim that for every $f\in\mL_T^\infty\bC^{-1,*}_{m,\eps}$,
$$
f*\phi_n\to f\ \hbox{in }\mL_T^\infty\bC^{-1,**}_{m,\eps}.
$$
Indeed, for the high frequencies, uniformly in $n$,
$$
\sup_{j\geq N}w_{m,\eps,0}(j)2^{-j}
\|\cR_j(f*\phi_n-f)\|_{\mL_T^\infty L^\infty}
\lesssim
\ff1{\ell_{m+2}(N)}
\|f\|_{\mL_T^\infty\bC^{-1,*}_{m,\eps}},
$$
which tends to $0$ as $N\to\infty$. For each fixed $N$, the finite number of dyadic blocks $j<N$ converges directly by mollification and the Bernstein inequality. Applying this to $f=b$ and $f=b\cdot\nn u^\l$, we obtain
$$
\|b-b*\phi_n\|_{\mL_T^\infty\bC^{-1,**}_{m,\eps}}\to0,
\ 
\|\phi_n*(b\cdot\nn u^\l)-b\cdot\nn u^\l\|_{\mL_T^\infty\bC^{-1,**}_{m,\eps}}\to0.
$$
On the other hand, the dyadic resolvent estimate gives
$$
w_{m,\eps,0}(j)2^j
\|\cR_j u^\l\|_{\mL_T^\infty L^\infty}
\lesssim
\ff{2^{2j}}{\l+2^{2j}}
\ff1{\ell_{m+2}(j)}
\|b+b\cdot\nn u^\l\|_{\mL_T^\infty\bC^{-1,*}_{m,\eps}}.
$$
Consequently the high frequency tail of $u^\l$ vanishes in $\mL_T^\infty\bC^{1,**}_{m,\eps}$, and the finite dyadic blocks converge by mollification. Thus
$$
\|u_n^\l-u^\l\|_{\mL_T^\infty\bC^{1,**}_{m,\eps}}\to0.
$$
By the paraproduct estimate in Lemma \ref{driftterm}, we obtain
$$
\|b\cdot\nn(u_n^\l-u^\l)\|_{\mL_T^\infty\bC^{-1,*}_{m,\eps}}
\lesssim
\|b\|_{\mL_T^\infty\bC^{-1,*}_{m,\eps}}
\|u_n^\l-u^\l\|_{\mL_T^\infty\bC^{1,**}_{m,\eps}}
\to0.
$$
Combining the above estimates, we get
$$
\|r_n\|_{\mL_T^\infty\bC^{-1,**}_{m,\eps}}\to0.
$$
Moreover, each $r_n$ is approximable in $\mL_T^\infty\bC^{-1,**}_{m,\eps}$ by the same mollification argument. Therefore the extended Krylov estimate for additive functionals yields, for every $p>2$,
$$
\sup_{t\in[0,T]}\|A_t^{r_n}\|_{L^p(\Om)}
\lesssim
\om_{**}(T)
\|r_n\|_{\mL_T^\infty\bC^{-1,**}_{m,\eps}}
\to0.
$$

Since $u_n^\l\to u^\l$ and $\nn u_n^\l\to\nn u^\l$ uniformly on $[0,T]\times\mR^d$, we can pass to the limit and get
\begin{align}
    \Phi(t,X_t)
    =\Phi(0,X_0)+\l\int_0^t u^\l(s,X_s)\dif s
    +\sqrt2\int_0^t\nn\Phi(s,X_s)\dif W_s.\no
\end{align}
Indeed, the convergence of the stochastic integral follows from the Burkholder-Davis-Gundy inequality:
\begin{align}
    \mE\bigg|\int_0^t(\nn\Phi_n-\nn\Phi)(s,X_s)\dif W_s\bigg|^2
    \leq T\|\nn\Phi_n-\nn\Phi\|_\infty^2\to0.\no
\end{align}
Since $X_t=\Phi^{-1}(t,Y_t)$, this is exactly the transformed SDE with coefficients $\wt b$ and $\wt\s$.

It remains to verify the stated properties of the coefficients. Boundedness is immediate from the boundedness of $u^\l$ and $\nn\Phi$. For the continuity, first note that if $0<|x-x'|=r\leq1/2$, then by the definition of $\bC^{1,**}_{m,\eps}$, using $J_{r^{-2}}$ as in \eqref{Jl-def},
\begin{align}
    |\nn u^\l(t,x)-\nn u^\l(t,x')|
    &\lesssim
    r\sum_{j\leq J_{r^{-2}}}2^j\|\cR_j\nn u^\l(t)\|_\infty
    +\sum_{j>J_{r^{-2}}}\|\cR_j\nn u^\l(t)\|_\infty\no\\
    &\lesssim
    \ff1{\ell_{m+1}^{\eps}(J_{r^{-2}})}
    \|u^\l\|_{\mL_T^\infty\bC^{1,**}_{m,\eps}}\lesssim
    \ff1{\ell_{m+2}^{\eps}(r^{-1})}
    \|u^\l\|_{\mL_T^\infty\bC^{1,**}_{m,\eps}}.\no
\end{align}
For $r>1/2$, the same bound follows from $\|\nn u^\l\|_\infty<\infty$, after enlarging the implicit constant. Now set
\begin{align}
    x:=\Phi^{-1}(t,y),\  x':=\Phi^{-1}(t,y').\no
\end{align}
Since $|x-x'|\leq4|y-y'|$, the slow variation of the iterated logarithms gives
\begin{align}
    |\wt b(t,y)-\wt b(t,y')|
    &\leq \l\|\nn u^\l\|_\infty |x-x'|
    \lesssim |y-y'|,\no\\
    |\wt\s(t,y)-\wt\s(t,y')|
    &\lesssim |\nn u^\l(t,x)-\nn u^\l(t,x')|\lesssim
    1\wedge\ff1{\ell_{m+2}^{\eps}(|y-y'|^{-1})}.\no
\end{align}
Finally, since $\|\nn u^\l\|_\infty\leq1/2$,
\begin{align}
    \ff12|\xi|\leq |\nn\Phi(t,x)\xi|\leq\ff32|\xi|,
    \  \xi\in\mR^d.\no
\end{align}
The proof is complete.
\epf

\bt\label{main-uniqueness}
Assume that $\eps>0$, $m\in\mN$, $b\in\mL_T^\infty\bC^{-1,*}_{m,\eps}$ and $\div b=0$. Weak uniqueness holds for SDE \eqref{SDE} in the class of weak solutions satisfying Definition \ref{krylov-class-def}. 
\et

\bpf
Let $(\mathfrak F,X,W)$ and $(\bar{\mathfrak F},\bar X,\bar W)$ be two weak solutions to \eqref{SDE} with the same initial law $\mu$, both satisfying Definition \ref{krylov-class-def}. By Lemma \ref{zvonkin-transform},
\begin{align}
    Y_t:=\Phi(t,X_t),\  \bar Y_t:=\Phi(t,\bar X_t)\no
\end{align}
solve the transformed SDE with the same coefficients $\wt b,\wt\s$ and the same initial law. By Lemma \ref{zvonkin-transform}, $\wt b$ is bounded Lipschitz, while $\wt\s$ is bounded, uniformly continuous and uniformly non-degenerate. Hence, by the classical weak uniqueness result for uniformly elliptic diffusions with bounded Lipschitz drift and uniformly continuous diffusion coefficient (cf. \cite[Theorem 7.2.1]{SV06}), the transformed SDE admits weak uniqueness. Therefore $\Law(Y)=\Law(\bar Y)$ on $C([0,T];\mR^d)$. Since $x\mapsto\Phi(t,x)$ is a diffeomorphism for each $t$, this implies $\Law(X)=\Law(\bar X)$ on $C([0,T];\mR^d)$. The proof is complete.
\epf

\subsection{Two-sided Gaussian bounds}

We conclude the proof of the main theorem by establishing the two-sided Gaussian bounds for the solutions starting from deterministic initial points.

\bp\label{gaussian-bounds}
Assume that $\eps>0$, $m\in\mN$, $b\in\mL_T^\infty\bC^{-1,*}_{m,\eps}$ and $\div b=0$.  For each $x\in\mR^d$, let $X^x$ denote the weak solution to \eqref{SDE} with $X^x_0=x$. Then, for every $t\in(0,T]$, the law of $X^x_t$ admits a density $\rho_t(x,\cdot)$ with respect to the Lebesgue measure, namely,
\begin{align}
    \bP(X^x_t\in A)=\int_A \rho_t(x,y)\dif y,\ A\in\cB(\mR^d).
    \no
\end{align}
The density $\rho_t(x,y)$ satisfies the following two-sided Gaussian estimate: there exist constants $C_0,C_1,\gamma_0,\gamma_1>0$, depending only on $T,d,m,\eps$ and $\|b\|_{\mL_T^\infty\bC^{-1,*}_{m,\eps}}$, such that for all $t\in(0,T]$, $x\in\mR^d$ and a.e. $y\in\mR^d$,
\begin{align}
    C_0t^{-\frac d2}\exp\Big(-\gamma_0\frac{|x-y|^2}{t}\Big)
    \leq \rho_t(x,y)
    \leq
    C_1t^{-\frac d2}\exp\Big(-\gamma_1\frac{|x-y|^2}{t}\Big).
    \no
\end{align}
\ep

\bpf
For simplicity, in this paragraph $\rho_n(t,x,y)$ and $\rho_t(x,y)$ denote the transition densities from time $0$ to time $t$.

Define, for $1\leq i,j\leq d$,
$$
B_{ij}:=(1-\Delta)^{-1}(\p_i b_j-\p_j b_i),\ b^0:=(1-\Delta)^{-1}b .
$$
Then $B_{ij}=-B_{ji}$. Moreover, since $\div b=0$,
$$
\sum_{j=1}^d\p_jB_{ij}=(1-\Delta)^{-1}\Big(\p_i\sum_{j=1}^d\p_jb_j-\Delta b_i\Big)=-(1-\Delta)^{-1}\Delta b_i,
$$
and hence
$$
b_i=\sum_{j=1}^d\p_jB_{ij}+b_i^0\ \hbox{in }\sS'(\mR^d).
$$
By the Bernstein inequality, for $j\geq1$,
$$
\|\cR_jB\|_{\mL_T^\infty L^\infty}\lesssim2^{-j}\|\cR_jb\|_{\mL_T^\infty L^\infty}\lesssim\ff{\|b\|_{\mL_T^\infty\bC^{-1,*}_{m,\eps}}}{w_{m,\eps,1}(j)}.
$$
Since $\sum_{j\geq0}w_{m,\eps,1}(j)^{-1}<\infty$, we obtain
$$
\|B\|_{\mL_T^\infty L^\infty}\lesssim\|b\|_{\mL_T^\infty\bC^{-1,*}_{m,\eps}}.
$$
Similarly,
$$
\|b^0\|_{\mL_T^\infty L^\infty}\lesssim\sum_{j\geq0}\|\cR_jb^0\|_{\mL_T^\infty L^\infty}\lesssim\sum_{j\geq0}2^{-2j}\|\cR_jb\|_{\mL_T^\infty L^\infty}\lesssim\sum_{j\geq0}\ff{2^{-j}\|b\|_{\mL_T^\infty\bC^{-1,*}_{m,\eps}}}{w_{m,\eps,1}(j)}\lesssim\|b\|_{\mL_T^\infty\bC^{-1,*}_{m,\eps}}.
$$
Thus
$$
\|B\|_{\mL_T^\infty L^\infty}+\|b^0\|_{\mL_T^\infty L^\infty}\lesssim\|b\|_{\mL_T^\infty\bC^{-1,*}_{m,\eps}}.
$$
Let $b_n=b*\phi_n$, $B_n=B*\phi_n$ and $b_n^0=b^0*\phi_n$. Then
$$
b_{n,i}=\sum_{j=1}^d\p_jB_{n,ij}+b_{n,i}^0.
$$
Moreover,
$$
\|B_n\|_{\mL_T^\infty L^\infty}+\|b_n^0\|_{\mL_T^\infty L^\infty}\leq\|B\|_{\mL_T^\infty L^\infty}+\|b^0\|_{\mL_T^\infty L^\infty}\lesssim\|b\|_{\mL_T^\infty\bC^{-1,*}_{m,\eps}},
$$
uniformly in $n$.

For every smooth $\varphi$,
$$
(\Delta-b_n\cdot\nn)\varphi
=\div\big((I-B_n^T)\nn\varphi\big)-b_n^0\cdot\nn\varphi.
$$
Indeed,
$$
\begin{aligned}
\div\big((I-B_n^T)\nn\varphi\big)
&=\Delta\varphi-\sum_{i,j=1}^d\p_i(B_{n,ji}\p_j\varphi)\\
&=\Delta\varphi-\sum_{j=1}^d\Big(\sum_{i=1}^d\p_iB_{n,ji}\Big)\p_j\varphi\\
&=\Delta\varphi-(b_n-b_n^0)\cdot\nn\varphi,
\end{aligned}
$$
where $\sum_{i,j}B_{n,ji}\p_{ij}\varphi=0$ by antisymmetry. Moreover,
$$
\xi\cdot (I-B_n^T)(t,x)\xi=|\xi|^2,\quad
|(I-B_n^T)(t,x)|
\leq1+\|B_n\|_{\mL_T^\infty L^\infty}
\lesssim1+\|b\|_{\mL_T^\infty\bC^{-1,*}_{m,\eps}}.
$$
Consequently, $\Delta-b_n\cdot\nn$ admits a divergence-form representation with principal coefficient matrix $I-B_n^T$ and first-order coefficient $-b_n^0$. Note that $I-B_n^T$ is uniformly elliptic and that the $L^\infty$-norms of $I-B_n^T$ and $b_n^0$ are bounded uniformly in $n$ by a constant depending only on $d,m,\eps$ and $\|b\|_{\mL_T^\infty\bC^{-1,*}_{m,\eps}}$. These uniform bounds allow us to apply Aronson's theorem \cite[Theorem 10(ii)]{Aro68} to $\partial_t-\operatorname{div}((I-B_n^T)\nabla)+b_n^0\cdot\nabla$. Hence its fundamental solution satisfies two-sided Gaussian bounds with constants independent of $n$. More precisely, for all $t\in(0,T]$ and $x,y\in\mR^d$,
$$
C_0t^{-\frac d2}\exp\Big(-\gamma_0\ff{|x-y|^2}{t}\Big)\leq\rho_n(t,x,y)\leq C_1t^{-\frac d2}\exp\Big(-\gamma_1\ff{|x-y|^2}{t}\Big),
$$
where the constants $C_0,C_1,\gamma_0,\gamma_1>0$ depend only on $T,d,m,\eps$ and $\|b\|_{\mL_T^\infty\bC^{-1,*}_{m,\eps}}$, and are independent of $n$.

We now pass to the limit. Fix $x\in\mR^d$. By the compactness argument in Theorem \ref{main-existence}, every subsequence of $\{X^{n,x}\}_{n\geq1}$ is tight in $C([0,T];\mR^d)$, and every subsequential limit is a weak solution to \eqref{SDE} with initial value $x$. Moreover, by Proposition \ref{constructed-solution-krylov}, such a limit satisfies Definition \ref{krylov-class-def}. The weak uniqueness already proved therefore identifies all subsequential limits with the law of $X^x$. Hence
$$
X^{n,x}\Rightarrow X^x\ \hbox{in }C([0,T];\mR^d).
$$
Thus, for every $t\in(0,T]$ and every nonnegative $\psi\in C_c(\mR^d)$,
$$
\int_{\mR^d}\psi(y)\rho_n(t,x,y)\dif y=\mE\psi(X_t^{n,x})\to\mE\psi(X_t^x).
$$
Set
$$
G_i(t,x,y):=C_it^{-\frac d2}\exp\Big(-\gamma_i\ff{|x-y|^2}{t}\Big),\ i=0,1.
$$
Passing to the limit in the preceding two-sided bound gives
$$
\int_{\mR^d}\psi(y)G_0(t,x,y)\dif y\leq\mE\psi(X_t^x)\leq\int_{\mR^d}\psi(y)G_1(t,x,y)\dif y.
$$
Since this holds for all nonnegative $\psi\in C_c(\mR^d)$, the usual regularity argument yields the measure inequality
$$
G_0(t,x,y)\dif y\leq\Law(X_t^x)(\dif y)\leq G_1(t,x,y)\dif y.
$$
In particular, $\Law(X_t^x)$ is absolutely continuous with respect to the Lebesgue measure. Let $\rho_t(x,\cdot)$ be its density. The preceding measure inequality implies
$$
G_0(t,x,y)\leq\rho_t(x,y)\leq G_1(t,x,y)
$$
for a.e. $y\in\mR^d$. This proves the desired Gaussian estimate.
\epf


\begin{bibdiv}
\begin{biblist}
\bibitem{Aro68}
D. G. Aronson, Non-negative solutions of linear parabolic equations, \emph{Ann. Scuola Norm. Sup. Pisa Cl. Sci.}, vol. 22, no. 4, pp. 607--694, 1968. 
\bibitem{BCD11}
H. Bahouri, J.-Y. Chemin and R. Danchin, \emph{Fourier analysis and nonlinear partial differential equations}, Grundlehren der mathematischen Wissenschaften [Fundamental Principles of Mathematical Sciences], Heidelberg, Springer, 2011.
\bibitem{BC01}
R. F. Bass and Z.-Q. Chen, Stochastic differential equations for Dirichlet processes, \emph{Probab. Theory Relat. Fields}, vol. 121, no. 3, pp. 422-446, 2001. \url{https://doi.org/10.1007/s004400100151}
\bibitem{BC03}
R. F. Bass and Z.-Q. Chen, Brownian motion with singular drift, \emph{Ann. Probab.}, vol. 31, no. 2, pp. 791--817, 2003. \url{https://doi.org/10.1214/aop/1048516536}
\bibitem{BFGM19}
L. Beck, F. Flandoli, M. Gubinelli and M. Maurelli, Stochastic ODEs and stochastic linear PDEs with critical drift: regularity, duality and uniqueness, \emph{Electron. J. Probab.}, vol. 24, 2019, Art. no. 136. \url{https://doi.org/10.1214/19-EJP379}
\bibitem{CC18}
G. Cannizzaro and K. Chouk, Multidimensional SDEs with singular drift and universal construction of the polymer measure with white noise potential, \emph{Ann. Probab.}, vol. 46, no. 3, pp. 1710--1763, 2018. \url{https://doi.org/10.1214/17-AOP1213}
\bibitem{CDFM15}
F. Colombini, D. Del Santo, F. Fanelli, and G. Métivier, The well-posedness issue in Sobolev spaces for hyperbolic systems with Zygmund-type coefficients, \emph{Comm. Partial Differential Equations},  vol. 40, no. 11, pp. 2082-2121, 2015. \url{https://doi.org/10.1080/03605302.2015.1082107}
\bibitem{DD16}
F. Delarue and R. Diel, Rough paths and $1d$ SDE with a time dependent distributional drift: application to polymers, \emph{Probab. Theory Relat. Fields}, vol. 165, no. 1--2, pp. 1--63, 2016. \url{https://doi.org/10.1007/s00440-015-0626-8}
\bibitem{Fol81}
H. Föllmer, Calcul d'Itô sans probabilités, \emph{Séminaire de probabilités de Strasbourg}, vol. 15, pp. 143-150, 1981. \url{https://doi.org/10.1007/BFb0088364}
\bibitem{FIR17}
F. Flandoli, E. Issoglio and F. Russo, Multidimensional stochastic differential equations with distributional drift, \emph{Trans. Amer. Math. Soc.}, vol. 369, no. 3, pp. 1665--1688, 2017. \url{https://doi.org/10.1090/tran/6729}
\bibitem{GP24}
L. Gr\"afner and N. Perkowski, Weak well-posedness of energy solutions to singular SDEs with supercritical distributional drift, arXiv:2407.09046, 2024.
\bibitem{HWZ20}
Z. Hao, M. Wu and X. Zhang, Schauder estimates for nonlocal kinetic equations and applications, \emph{J. Math. Pures Appl. (9)}, vol. 140, pp. 139-184, 2020. \url{https://doi.org/10.1016/j.matpur.2020.06.003}
\bibitem{HZ25}
Z. Hao and X. Zhang, SDEs with supercritical distributional drifts, \emph{Comm. Math. Phys.}, vol. 406, no. 10, 2025, Art. no. 250. \url{https://doi.org/10.1007/s00220-025-05430-2}
\bibitem{KP22}
H. Kremp and N. Perkowski, Multidimensional SDE with distributional drift and L\'evy noise, \emph{Bernoulli}, vol. 28, no. 3, pp. 1757--1783, 2022. \url{https://doi.org/10.3150/21-BEJ1394}
\bibitem{Kry21}
N. V. Krylov, On diffusion processes with drift in $L_d$, \emph{Probab. Theory Relat. Fields}, vol. 179, no. 1--2, pp. 165--199, 2021. \url{https://doi.org/10.1007/s00440-020-01007-3}
\bibitem{KR05}
N. V. Krylov and M. R\"ockner, Strong solutions of stochastic equations with singular time dependent drift, \emph{Probab. Theory Relat. Fields}, vol. 131, no. 2, pp. 154--196, 2005. \url{https://doi.org/10.1007/s00440-004-0361-z}
\bibitem{Le20}
K. Lê, A stochastic sewing lemma and applications, \emph{Electron. J. Probab.}, vol. 25, pp. 1–55, 2020. \url{https://doi.org/10.1214/20-EJP442}
\bibitem{MP24}
T. Matsuda and N. Perkowski, An extension of the stochastic sewing lemma and applications to fractional stochastic calculus, \emph{Forum Math. Sigma}, vol. 12, 2024, Art. no. e52. \url{https://doi.org/10.1017/fms.2024.32}
\bibitem{RZ23}
M. R\"ockner and G. Zhao, SDEs with critical time dependent drifts: weak solutions, \emph{Bernoulli}, vol. 29, no. 1, pp. 757--784, 2023. \url{https://doi.org/10.3150/22-BEJ1478}
\bibitem{SSSZ12}
G. Seregin, L. Silvestre, V. \v{S}ver\'ak and A. Zlato\v{s}, On divergence-free drifts, \emph{J. Differential Equations}, vol. 252, no. 1, pp. 505--540, 2012. \url{https://doi.org/10.1016/j.jde.2011.08.039}
\bibitem{SV06}
D. Stroock and S. Varadhan, \emph{Multidimensional diffusion processes}, Classics in Mathematics. Berlin, Springer-Verlag, 2006.
\bibitem{Ver79}
A. Ju. Veretennikov, On the strong solutions of stochastic differential equations, \emph{Theory Probab. Appl.}, vol. 24, no. 2, pp. 354--366, 1979. \url{https://doi.org/10.1137/1124039}
\bibitem{Zha05}
X. Zhang, Strong solutions of SDEs with singular drift and Sobolev diffusion coefficients, \emph{Stoch. Process. Appl.}, vol. 115, no. 11, pp. 1805--1818, 2005. \url{https://doi.org/10.1016/j.spa.2005.06.003}
\bibitem{Zha11}
X. Zhang, Stochastic homeomorphism flows of SDEs with singular drifts and Sobolev diffusion coefficients, \emph{Electron. J. Probab.}, vol. 16, pp. 1096--1116, 2011. \url{https://doi.org/10.1214/EJP.v16-887}
\bibitem{Zha16}
X. Zhang, Stochastic differential equations with Sobolev diffusion and singular drift and applications, \emph{Ann. Appl. Probab.}, vol. 26, no. 5, pp. 2697--2732, 2016. \url{https://doi.org/10.1214/15-AAP1159}
\bibitem{Zha24}
X. Zhang, Maximum principle for non-uniformly parabolic equations and applications, \emph{Ann. Sc. Norm. Super. Pisa Cl. Sci.}, vol. 25, no. 3, pp. 1267--1308, 2024. \url{https://doi.org/10.2422/2036-2145.202105_052}
\bibitem{ZZ18}
X. Zhang and G. Zhao, Singular Brownian diffusion processes, \emph{Commun. Math. Stat.}, vol. 6, no. 4, pp. 533--581, 2018. \url{https://doi.org/10.1007/s40304-018-0164-7}
\bibitem{ZZ21}
X. Zhang and G. Zhao, Stochastic Lagrangian path for Leray's solutions of 3D Navier--Stokes equations, \emph{Comm. Math. Phys.}, vol. 381, no. 2, pp. 491--525, 2021. \url{https://doi.org/10.1007/s00220-020-03888-w}
\end{biblist}
\end{bibdiv}
\end{document}